\documentclass[11pt]{amsart}
\usepackage[T1]{fontenc}
\usepackage{lmodern,mathtools,amssymb}
\usepackage[margin=1.12in]{geometry}
\usepackage{xcolor}
\definecolor{definitionpurple}{RGB}{94,38,125}
\DeclareRobustCommand{\defn}[1]{\textcolor{definitionpurple}{\textit{#1}}}
\usepackage[hidelinks]{hyperref}
\DeclareRobustCommand{\thmbound}[2]{f_{\hyperref[#1]{\ref*{#1}}}\!\left(#2\right)}
\newtheorem{theorem}{Theorem}[section]
\newtheorem{lemma}[theorem]{Lemma}
\newtheorem{proposition}[theorem]{Proposition}
\title{The $k$-disjoint paths problem for group-labelled graphs}
\author{Tony Huynh}
\address{Discrete Mathematics Group, Institute for Basic Science (IBS), Daejeon, Republic of Korea}
\thanks{The author is currently supported by the Institute for Basic Science (IBS-R029-C1) and was supported by an NSERC grant while he was a PhD student.}
\hypersetup{pdftitle={The k-disjoint paths problem for group-labelled graphs},pdfauthor={Tony Huynh}}
\date{}
\subjclass[2020]{05C38, 05C83, 05C85, 05C10}
\keywords{disjoint paths, gain graph, group-labelled minor, graph embedding}
\begin{document}
\begin{abstract}
A $\Gamma$-labelled graph is an oriented graph with edge labels in a group $\Gamma$.
Fix a finite abelian group $\Gamma$ and a positive integer $k$. We give a polynomial-time
algorithm that finds $k$ vertex-disjoint paths joining prescribed pairs with a prescribed
gain on each path, or determines that none exist. The result previously appeared only in the author's
2009 PhD thesis. This proof uses a reduction to unbreakable graphs and recent
structure theorems for graphs drawn on a surface, and was found with the help of AI.
\end{abstract}
\maketitle
\section{Introduction}
Let $\Gamma$ be a finite abelian group written additively. An \defn{oriented graph}
is a finite graph in which each edge has a designated \defn{tail} and \defn{head}. We allow
loops and parallel edges. A \defn{$\Gamma$-labelled graph} $G$ is an oriented graph with
a \defn{label} $\gamma_G(e)\in\Gamma$ for each edge $e$. Write $\operatorname{tail}_G(e)$
and $\operatorname{head}_G(e)$ for its tail and head. We write $\vec G$ for the directed
graph obtained by forgetting the labels and $\widetilde G$ for the \defn{underlying undirected
graph}, obtained by also forgetting the orientations. Connectivity, embeddings and treewidth of $G$ refer to $\widetilde G$.

A \defn{walk} in $G$ is a sequence $W\coloneqq(v_0,e_1,v_1,\ldots,e_\ell,v_\ell)$
in which $e_i$ joins $v_{i-1}$ to $v_i$ in $\widetilde G$ for $1\le i\le\ell$.
Vertices and edges may recur. At a loop the walk also specifies the end incidence at
which it arrives. Write $\operatorname{tail}_G(W)\coloneqq v_0$ and
$\operatorname{head}_G(W)\coloneqq v_\ell$. For an end incidence $\iota$ of an edge $e$,
put $\gamma_G(e,\iota)\coloneqq\gamma_G(e)$ at the head and
$\gamma_G(e,\iota)\coloneqq-\gamma_G(e)$ at the tail. For a non-loop we may write its
incident vertex in place of $\iota$. If $\iota_i$ is the arrival incidence of $e_i$,
the \defn{gain} of $W$ is
\[
 \gamma_G(W)\coloneqq\sum_{i=1}^{\ell}\gamma_G(e_i,\iota_i).
\]
A \defn{path} is a walk with no repeated vertex. A \defn{cycle} is a nonempty closed walk with distinct edges and
no repeated vertex except its ends. A loop and a pair of parallel edges can thus give
cycles of lengths one and two, respectively. \defn{Flipping} an edge interchanges its
head and tail and negates its label. Corresponding traversals have the same gain.

For an integer $h\ge0$, put $[h]\coloneqq\{1,\ldots,h\}$, with $[0]\coloneqq\varnothing$.
For an integer $k\ge0$, a \defn{$k$-pattern} in $G$ is a set
$\Pi\coloneqq \{(s_i,t_i,\gamma_i):i\in[h]\}$ with $0\le h\le k$, where the $2h$ vertices
$s_i,t_i$ are distinct vertices of $G$ and each $\gamma_i$ belongs to $\Gamma$.
Write $V(\Pi)\coloneqq \{s_i,t_i:i\in[h]\}$ for its set of \defn{terminals}.
We call each triple an \defn{ordinary demand}. A \defn{$\Pi$-linkage} in $G$ is a
family of pairwise vertex-disjoint paths $P_i$, $i\in[h]$, such that
\[
  \operatorname{tail}_G(P_i)=s_i,\qquad \operatorname{head}_G(P_i)=t_i,\qquad
  \gamma_G(P_i)=\gamma_i .
\]
The \defn{$k$-disjoint paths problem} asks for such a linkage.

A \defn{shift} by a \defn{potential} $p:V(G)\to\Gamma$ adds
$p(\operatorname{head}_G(e))-p(\operatorname{tail}_G(e))$ to the label of each edge $e$.
The \defn{size} of $G$ is $n\coloneqq |V(G)|+|E(G)|$. Throughout the paper, the subscript
on a bound $f$ gives the number of the result supplying it.

\begin{theorem}\label{thm:main}
Let $\Gamma$ be a fixed finite abelian group and let $k$ be a fixed positive integer.
There is an integer $\thmbound{thm:main}{\Gamma,k}$ such that, given a $\Gamma$-labelled
graph $G$ of size $n$ and a $k$-pattern $\Pi$ in $G$, a deterministic algorithm runs in time
$O(n^{\thmbound{thm:main}{\Gamma,k}})$ and finds a $\Pi$-linkage or determines that none exists.
\end{theorem}

As observed in the thesis~\cite[Theorems~1.1.10 and~1.4.9]{HuynhThesis},
Theorem~\ref{thm:main} also gives polynomial-time algorithms for testing whether $G$
contains a fixed $\Gamma$-labelled graph as a minor or as a subdivision.
Here a \defn{labelled minor} is obtained by shifts, deletions and contractions of
zero-labelled edges. A \defn{labelled subdivision} replaces each edge by a path
(a cycle for a loop) of the same gain in its given direction, using new internal
vertices private to that edge. Geelen and Gerards posed the minor-testing
problem~\cite[Conjecture~1.2]{GeelenGerards}.

\subsection*{Related Work}
The ordinary disjoint paths theorem of Robertson and Seymour~\cite{RSXIII} is the case of the
trivial group. In the unoriented model the gain of a path is the sum of its edge
labels. Reversing a path preserves that sum and negates its gain in the oriented model.
The two conventions therefore agree when every group element has order at most two.
For the unoriented model, Liu and Yoo~\cite[Theorem~1.3]{LiuYoo} give an
$O(|V(G)|^8+|E(G)|)$-time algorithm for fixed $\Gamma$ and $k$. To prescribe gain
$\gamma_i$, forbid the set $\Gamma\setminus\{\gamma_i\}$ on path $i$.

When the input contained a large grid as a minor, the thesis used the grid to
certify a linkage or find a vertex whose deletion preserved the answer. After deleting
such a vertex, the algorithm restarted on the smaller graph~\cite[Section~1.5 and Chapter~10]{HuynhThesis}.
The 2018 model-checking theorem of Lokshtanov, Ramanujan, Saurabh and Zehavi~\cite{LRSZ}
reduces the decision problem to graphs in which every separation of a prescribed small
order has only boundedly many vertices on one side outside the separator.

Both proofs use the gain around a face to adjust the gain of a path. The difference
is how we find enough places to make these adjustments. We use the 2025 finite-index
structure theorem of Paul, Protopapas, Thilikos and Wiederrecht~\cite{PPTWindex}
in the surface forms of Bla\v{z}ej, Pilipczuk and Protopapas~\cite{BPPcoarseMenger}
from 2026. We record the gains of incident faces, up to sign, as colours on pieces
of the drawing. The resulting structure guarantees that each nonzero facial gain
either occurs in many disjoint disks where we can adjust a path, or has all its
occurrences confined to a bounded collection of disks. The gains available for
these adjustments determine the subgroup used in the surface algorithm.

\subsection*{Outline of the paper}

Section~\ref{sec:linkage-basics} gives the linkage reductions used throughout. It
encodes prescribed gains in monadic second-order logic, handles graphs of bounded
treewidth, and records how a linkage can enter and leave a specified vertex set.
Section~\ref{sec:flatness} tests where a labelling on an embedded graph extends
across the surrounding surface. Section~\ref{sec:journal-attachments} bounds the
number of passages a linkage needs between a surface drawing and the graphs
attached along its boundary.

Section~\ref{sec:embedded-reductions} reduces an embedded instance to
3-connected pieces and proves the assertion above about nonzero facial gains.
Section~\ref{sec:surface-routing} builds a finite embedded graph that carries
all the required relative homology classes. A path can follow either side of
each assigned face, and changing sides adds its facial gain.
Section~\ref{sec:surface-algorithm} combines these observations into the
bounded-genus algorithm. If the available facial gains generate the current
subgroup, we correct the paths independently. Otherwise we reduce to a smaller
subgroup while keeping the original target gains.

Section~\ref{sec:journal-global} reduces arbitrary graphs to bounded genus.
A large complete minor supplies prescribed gains on a small separator. In the
excluded-minor case Appendix~\ref{app:construct} constructs the
surface drawing and the graphs attached along its boundary. It uses the
near-embedding theorem of Kawarabayashi, Thomas and Wollan~\cite{KTW} and the
refinement of Diestel, Kawarabayashi, M\"uller and Wollan~\cite{DKMW}, in the
algorithmic forms stated by Fiorini, Joret, Weltge and
Yuditsky~\cite[Theorems~15 and~19]{FJWY}. The section ends by recovering the
paths from the decision algorithm.

The algorithms are nonuniform in $\Gamma$ and $k$. The reduction of
Lokshtanov, Ramanujan, Saurabh and Zehavi~\cite{LRSZ} in
Section~\ref{sec:linkage-basics} supplies a separation bound for each fixed
formula. Section~\ref{sec:surface-routing} uses a bound on the size of a
finite embedded graph. We also use the finite list of excluded minors for
each surface, whose existence follows from Robertson and Seymour~\cite{RSXX}.
These constants enter the algorithm for the fixed parameters. The proof
does not compute them as functions of the parameters.

\section{Linkages and separations}\label{sec:linkage-basics}
An \defn{ordinary minor} is obtained by vertex deletions, edge deletions and edge contractions. A
\defn{subdivision} replaces edges by paths whose interiors are pairwise disjoint and avoid the
original vertices. A \defn{linear forest} is a graph whose components are paths, allowing isolated
vertices. A \defn{cutvertex} is a vertex whose deletion increases the number of connected components.
For linkage we discard loops before taking blocks. A \defn{block} on at least three vertices is a
maximal 2-connected subgraph, where \defn{2-connectivity} means that deleting any one vertex
leaves a connected graph. A pair of adjacent vertices lying in no such subgraph forms a
block with all edges between them. This includes a bridge and a bundle of parallel edges.
An isolated vertex is also a block. We handle blocks on at most two vertices directly.

A \defn{reserved pattern} in a $\Gamma$-labelled graph $G$ extends an ordinary pattern
by allowing \defn{singleton demands} $(x,x,0)$ for distinct vertices $x$ of $G$ outside its
ordinary demands. A \defn{realization} extends the linkage by the length-zero path at each such
$x$ and all its paths are vertex-disjoint. We say that a singleton \defn{reserves} its
vertex. The \defn{order} of a pattern is its number of demands. The vertices specified in the demands are the \defn{named vertices}. In the reductions below,
\defn{pattern} allows reservations and \defn{linkage} means a realization of the given
pattern. We may discard loops when searching for a linkage.

For a walk from $s$ to $t$, shifting by $p$ changes its gain by $p(t)-p(s)$, by telescoping.
We make this adjustment to every target when shifting. A shift can make the edges of any
forest zero by choosing a potential at one root of each tree and then at its children
recursively. Consequently all cycle gains in a connected graph belong to a subgroup $L$
precisely when a shift puts all edge labels in $L$. To see sufficiency of the cycle condition,
normalize a spanning tree. The \defn{fundamental cycle} of an edge outside the tree consists of that edge
and the tree path between its ends. After normalization its gain, up to sign, equals the
label of the edge.

A \defn{separation} of a graph with vertex set $V$ is a pair $(A,B)$ of vertex sets with
$A\cup B=V$ and no edge between $A\setminus B$ and $B\setminus A$. These two differences are
the \defn{strict sides}. The \defn{separator} is $A\cap B$, and the \defn{order} is $|A\cap B|$.
Equivalently one may distribute the edges between two subgraphs on these sets. A graph is
\defn{$(s,c)$-unbreakable} if every separation of order at most $c$ has a strict side with at
most $s$ vertices. A \defn{tree decomposition} consists of a tree with a vertex set of the
graph, called a \defn{bag}, at each node. Each graph vertex occurs on a nonempty connected subtree
and both ends of every edge occur in some common bag. The \defn{width} is the largest bag size minus
one. \defn{Treewidth} is the minimum width. A \defn{path decomposition} is one in which the tree is a
path. Its \defn{adhesion} is the largest size of an intersection of consecutive bags.

When we say that a problem \defn{reduces to} other problems, or that it suffices to solve them, we mean
that there is a polynomial-time algorithm for it which calls algorithms for the others
polynomially often. Such reductions are nested to a bounded depth in this paper.

The \defn{incidence encoding} of a loopless $\Gamma$-labelled graph subdivides every edge
once. It specifies the original vertices, partitions the subdivision vertices according to
the group labels, and marks the incidence joining each subdivision vertex to the tail of
its edge. Named pattern vertices are parameters. \defn{Monadic second-order logic} here allows
quantification over vertices, edges and sets of either. Unbreakability of a structure
refers to its underlying graph.

\begin{lemma}\label{lem:logic}
Let $\Gamma$ be a fixed finite abelian group and let $k$ be a fixed nonnegative integer.
On the incidence encodings of $\Gamma$-labelled graphs with patterns of order at most $k$,
finitely many fixed monadic second-order formulas express the existence of a linkage, one
for each demand type and target tuple. Let $\varphi$ be the conjunction of one of these
formulas with a fixed monadic second-order condition, and let $c$ be a fixed nonnegative
integer. There exists a nonnegative integer $s\coloneqq \thmbound{lem:logic}{\varphi,c}$ such that a polynomial-time decision
algorithm for $\varphi$ on $(s,c)$-unbreakable structures suffices to decide $\varphi$
on all structures in polynomial time.
\end{lemma}
\begin{proof}
The formula first checks the incidence encoding. Edges join original vertices to edge
vertices. Each edge vertex has two distinct neighbours, one marked tail incidence
and one label. Parallel edges have separate edge vertices. For each ordinary demand
$(s,t,\alpha)$ choose a connected edge set with degree one at $s,t$ and degree two
at its other incident vertices. Require these paths to be disjoint and to avoid all
reserved vertices. To check all gains, choose a single map $p:V(G)\to\Gamma$ by its
colour classes. On each selected edge $e$ directed from $u$ to $v$, require
$p(v)-p(u)=\gamma_G(e)$. At the ends of demand $i$, require $p(s_i)=0$ and
$p(t_i)=\gamma_i$. These equations telescope along each path in either
direction. Conversely, the partial sums on the disjoint paths give a consistent map,
which can be extended arbitrarily to unused vertices.

Apply the structure version of the model-checking theorem of Lokshtanov, Ramanujan,
Saurabh and Zehavi~\cite[Theorem~2]{LRSZ}, with free sets and individuals as parameters.
For fixed $\varphi,c$ it gives $s$ before the algorithm on unbreakable inputs is chosen.
The exponent of that algorithm may therefore depend on $s$ and can be increased to
exceed four as required. Finally, any separation of the decoded graph lifts to its
encoding with the same order and with all original vertices on their previous sides.
Place each edge vertex on a side containing both ends. Thus an unbreakable encoding
also has an unbreakable decoded graph with the same bounds.
\end{proof}

\begin{lemma}\label{prop:width}
Let $\Gamma$ be a fixed finite abelian group and let $w,k$ be fixed nonnegative integers.
Given a $\Gamma$-labelled graph $G$ of treewidth at most $w$ and a pattern $\Pi$ in $G$
of order at most $k$, one can decide in polynomial time whether $\Pi$ has a realization.
\end{lemma}
\begin{proof}
A width-$w$ decomposition of the input gives a decomposition of its encoding of width at most
$\max\{w,2\}$ by attaching a three-vertex bag for each subdivided edge to a bag containing
both ends. Use the fixed-width algorithm of Bodlaender~\cite{Bodlaender} and the bounded-treewidth
model checking theorems of Courcelle~\cite{Courcelle} and of Arnborg, Lagergren and
Seese~\cite{ALS} with Lemma~\ref{lem:logic}.
\end{proof}

For nonnegative integers $k,a$, put $\thmbound{near:apex}{k,a}\coloneqq k+a$.

\begin{lemma}\label{near:apex}
Let $\Gamma$ be a fixed finite abelian group and let $k,a$ be fixed nonnegative integers.
Let $G$ be a $\Gamma$-labelled graph, let $\Pi$ be a pattern in $G$ of order at most $k$,
and let $Z\subseteq V(G)$ have at most $a$ vertices. Deciding whether $\Pi$ has a
realization reduces to linkage decisions for patterns of order at most
$\thmbound{near:apex}{k,a}$ in the components of $G-Z$, with labels still in $\Gamma$.
\end{lemma}
\begin{proof}
For each ordinary path, guess its ordered visits to $Z$ and the incident edges used
at each visit. Use each visited vertex only once, omit reserved vertices, and use
one incidence at an endpoint. These choices cut the path into pieces in $G-Z$.
An edge joining consecutive visits leaves no intervening piece. Guess the gain of
each piece so that, together with the traversed edges, the gains sum to the target.

Test the pieces in each component together, keeping the original reservations there.
A one-vertex piece reserves its vertex and has gain zero. Reject inconsistent
incidences, conflicting terminals and pieces whose ends lie in different components.
Every linkage supplies a guess, and concatenation recovers a linkage from every
successful guess. Each visit adds at most one piece, so there are at most $k+|Z|$
demands and reservations in total. Only $2|Z|$ edge incidences need be guessed.
\end{proof}

\section{Gains on a punctured surface}\label{sec:flatness}
A \defn{surface} in this paper is a compact connected two-dimensional manifold, possibly
nonorientable. A surface is \defn{closed} if it has empty boundary. The \defn{Euler genus} of a surface with
$b$ boundary components and Euler characteristic $\chi$ is $2-b-\chi$, the Euler genus
obtained by capping each boundary component with a disk. The \defn{Euler genus of a graph} is the
least Euler genus of a closed surface in which it embeds. An \defn{embedding} is a drawing whose
edges meet only at common ends. A \defn{face} is a component of the complement. A \defn{cellular embedding}
is one in which every face is an open disk. A \defn{facial boundary walk} is the closed walk tracing the
boundary of a face in such an embedding. An edge can occur twice on this walk. Equivalently the closure of the face is the image of a closed
parameter disk under the map attaching its boundary walk to the graph. A \defn{noose} is a
simple closed curve meeting the drawing only in vertices. The \defn{representativity}, or \defn{face-width}, of a
drawing on a closed surface of positive Euler genus is the minimum number of intersections of
a noncontractible simple closed curve with the drawing. A \defn{corner} at a vertex is a local
sector between consecutive edge ends. A \defn{proper arc} in a surface with boundary is an
embedded interval that meets the boundary precisely at its two ends. Embeddings and punctures
can be specified by finite triangulations containing the graph as a subcomplex after
subdivision. Combinatorial embeddings give such triangulations of polynomial size by
subdividing regular neighbourhoods of the graph and triangulating the complementary surfaces.
In a noncellular drawing the incidence, orientability and genus data of those surfaces are
part of the embedding.

Let $G$ be a $\Gamma$-labelled graph drawn on a surface $S$. Let $Y$ be a finite set of
interior points of $S$ disjoint from $G$. To \defn{puncture} the surface at $Y$, remove a small open
disk around each point, with the closed disks pairwise disjoint and disjoint from $G$. A \defn{flat extension off
$Y$} assigns a gain to each path between points of the resulting surface, invariant under
homotopy with fixed endpoints, additive under concatenation and equal to the given gains on
edges of $G$. Constant paths have gain zero. We can equivalently specify transport on a
triangulation by assigning gains to its oriented edges so that reversal negates the gain and
the sum around each triangle is zero. In fact edge paths generate paths between triangulation
vertices up to homotopy, with the triangle and backtracking moves as relations. To include
other endpoints choose a reference path to each from a triangulation vertex. The gain of an
oriented closed curve is then invariant under free homotopy because $\Gamma$ is abelian. A
flat extension restricts to every subdrawing. For a subgroup $L\le\Gamma$, the \defn{quotient gain}
of a walk $W$ in $G$ is the coset $\gamma_G(W)+L$ in $\Gamma/L$. The \defn{exponent} $\exp(\Gamma/L)$ is the least positive
integer annihilating that quotient. In particular the trivial quotient has exponent one.

In a cellular drawing on a closed surface, the \defn{facial gain} of a face is the gain
of its boundary walk in a chosen direction. Reversing the walk negates this gain. A face
is \defn{defective} if its facial gain is nonzero. It is \defn{defective modulo} a subgroup
$L\le\Gamma$ if its facial gain has nonzero image in $\Gamma/L$.

\begin{lemma}\label{thm:find-flatness-punctures}
Let $\Gamma$ be a fixed finite abelian group. Given an embedding of a $\Gamma$-labelled
graph $G$ in a surface $\Sigma$, one can find in polynomial time a smallest set of
punctures off which the labelling extends flatly. If $\Sigma$ is closed and the embedding
is cellular, placing one puncture in each defective face gives such a set.
\end{lemma}
\begin{proof}
Triangulate the drawing, dividing the gain of each graph edge among its segments.
Any two such divisions differ by a shift fixing the original vertices. The unknown
edge labels belong to individual complementary regions, so consider one region at a time.

Glue copies of its triangles along just the sides crossed by a spanning tree of their dual graph.
This gives a disk. Labels on its boundary extend across its triangles precisely when
their oriented sum is zero. Indeed, peeling off leaf triangles determines the label
of each side leading to a parent, leaving just the boundary sum at the root.

The graph edges contribute a prescribed sum. A free edge on the boundary of
$\Sigma$ occurs once, so its label alone can cancel that sum. If there is no
such edge, each free edge occurs twice. Opposite signs cancel, while equal
signs contribute twice its label up to sign. The possible corrections are then
$2\Gamma\coloneqq\{2a:a\in\Gamma\}$ if some pair has equal signs, and $\{0\}$
otherwise. Testing whether a correction cancels the prescribed sum takes polynomial time.
A failed region needs a puncture. One suffices, since puncturing the root triangle
removes its equation and leaves every other step of the peeling unchanged.
For a cellular drawing on a closed surface, orient all triangles of a face by its
parameter disk. Every free pair then cancels, leaving exactly the facial gain.
\end{proof}

\section{Graphs attached along boundary arcs}\label{sec:journal-attachments}

Let $F$ be drawn in a compact connected surface $S$ with nonempty boundary. Edge interiors lie
in the interior of $S$. Take $r$ pairwise disjoint closed subarcs of the boundary, called
\defn{supports}, whose union contains every boundary vertex. Each is a proper subset of its
boundary circle and contains at least one drawing vertex. We allow
a support to be a single vertex and shrink a support with just one drawing vertex to that
vertex. Each boundary circle contains a drawing vertex. A nonsingleton support is
parametrized in one of its two linear orders. To a support with vertices $v_0,\ldots,v_m$ in
order attach a graph $H$ meeting $F$ precisely at these vertices. The graph $H$ is the
\defn{attachment}, and the ordered sequence of its vertices on the support is its \defn{society}.
Different attached graphs are vertex-disjoint. For $m\ge1$ a \defn{rooted presentation of
depth $b$} of $H$ consists of \defn{root sets} $X_0,\ldots,X_m$ and edge-disjoint \defn{pieces}
$K_0,\ldots,K_{m-1}$ covering $H$ with the following properties.
\begin{enumerate}
\item Each $X_i$ has $b$ elements and meets the society just at $v_i$. The graph $K_i$
contains $X_i\cup X_{i+1}$ and meets the society just at $v_i,v_{i+1}$.
\item For $i<j$, every vertex of $K_i\cap K_j$ belongs to each $X_h$ with $i<h\le j$.
\item There are $b$ vertex-disjoint paths in $K_i$ from $X_i$ to $X_{i+1}$. They may have length zero.
\end{enumerate}
On a singleton support the attachment consists of an edgeless root set $X_0$ of size $b$
containing its society vertex. Write $q$ for an upper bound on the depths. The local path
families concatenate into $b$ disjoint paths called \defn{rails}. Indeed, nonadjacent pieces meet
only at intervening roots, where every vertex is already used by both adjacent families. A
root persisting through several cuts has a trivial intervening path. Thus each rail meets each
$X_i$ in one vertex and encounters them in order, with repetitions allowed. A path used as an
attachment, with its vertices in path order along the support, has depth one.

In this section we use the subgraph convention for a separation. The two sides $(A,B)$ are
edge-disjoint subgraphs covering the graph.
Its order is $|V(A)\cap V(B)|$. We first record how to replace the part of a linkage between two
separations.

\begin{lemma}\label{lem:journal-pump}
Let $\Gamma$ be a finite abelian group and let $k,n$ be nonnegative integers. There is an
integer $m\coloneqq \thmbound{lem:journal-pump}{\Gamma,k,n}\ge2$ with the following property. Let $G$ be a $\Gamma$-labelled graph, and let
$(A_1,B_1),\ldots,(A_m,B_m)$ be separations of $G$ of order $n$, with
$A_i\subseteq A_{i+1}$ and $B_{i+1}\subseteq B_i$ for $1\le i<m$. Let
$L_1,\ldots,L_n$ be vertex-disjoint paths in $G$ from the first separator to the last.
Every realizable pattern in $G$ of order at most $k$ has a realization all of whose
edges in $B_i\cap A_j$ lie on $L_1,\ldots,L_n$ for some $1\le i<j\le m$.
\end{lemma}
\begin{proof}
Write $x_{ia}$ for the contact of $L_a$ with the $i$th separator. Disjointness and
cardinality give exactly one contact per path at each cut. They occur in order on each
$L_a$, with repetitions allowed. Shift all these paths to gain zero on every edge
and adjust the targets.

Fix a realization. At cut $i$, record whether each named vertex is strictly left, strictly
right or equal to $x_{ia}$ for a specified $a$. Restrict the realization to $A_i$ and mark
its named and separator vertices, including isolated marks. Every unmarked vertex incident
with an edge has degree two. Suppress these vertices. The record consists of the resulting
forest on the marks, with each separator vertex named by its index $a$, the gain of each of
its edges in path order, and for each mark the path of the realization using it, if any. This forest has at most $2k+n$
marks, so choosing $m$ larger than the number of possible records gives equal records
at some $i<j$.

Glue the prefix in $A_i$ to the suffix in $B_j$ by matching their separator contacts
rail by rail. As an abstract gluing this reproduces the original solution at cut $j$,
with the same names and gains. To realize it in $G$, use the segment of $L_a$ from
$x_{ia}$ to $x_{ja}$ whenever a path crosses between the two sides at this contact.
If both incidences are on one side, no segment is needed. Equal contacts already coincide.
All shared vertices of $A_i$ and $B_j$ are such equal contacts. The added segments
are pairwise disjoint, their interiors miss the retained sides, and the location
records ensure that no moving contact or segment interior is named. The segments add
zero gain. Undoing the shift gives the asserted realization.
\end{proof}

We will also use the following observation. On splitting vertices into copies, assign each old
edge end to a single copy. If the total increase in vertex count is $t$, any solution lifts to
a pattern of at most $k+3t$ demands and reservations for which every solution projects back.
Indeed, split its paths at the affected incidences, prescribe the resulting ordered pieces
with their actual gains, choose one copy of every original singleton, and reserve all unused
copies of split vertices. There are at most $k+t$ pieces and at most $2t$ additional copies to
reserve. Two used copies of the same vertex are necessarily the prescribed ends of two
successive pieces of the same original path. Identifying them thus makes exactly the intended
joins. At every other split vertex at most one copy is available. New solutions therefore
project to disjoint simple paths of the same gains.

For a pattern in a drawing $F$ with attachments as above, let $W$ consist of the society
vertices and all named vertices of $F$. A \defn{chord} of a realization is a positive-length
segment of one of its paths in $F$ with both ends in $W$ and no internal vertex in $W$.

\begin{lemma}\label{near:chords}
Fix a compact connected surface $S$ with nonempty boundary, a finite abelian group $\Gamma$,
and nonnegative integers $r,q,k$. There is an integer $M\coloneqq \thmbound{near:chords}{S,\Gamma,r,q,k}$ with the following property.
Let $G$ be a $\Gamma$-labelled graph presented as a drawing $F$ in $S$ with attachments
on at most $r$ supports as above, each of depth at most $q$. Every realizable pattern
$\Pi$ in $G$ of order at most $k$ has a realization using at most $M$ chords.
\end{lemma}
\begin{proof}
Put every named vertex of the drawing on a support. At each interior name remove a
small open disk in an incident face sector whose boundary touches the drawing only
at that name. Give it a singleton support of depth one. This leaves the chords
unchanged and adds at most $2k$ boundary circles and as many supports. Take the
maximum of the bounds below over these finitely many surface types, with support
and depth bounds $r+2k$ and $\max\{q,1\}$. Later cuts create names only on existing supports.

Choose a realization with the fewest chords. Retain its drawing edges and the society
vertices and suppress the internal vertices of each chord, summing gains.
The resulting drawing is a linear forest of proper chord edges and isolated vertices,
all on the boundary. Every solution in this graph expands in the original drawing.
Our chosen solution uses every chord edge, so a solution omitting even one would
contradict minimality.

First suppose $S$ is a disk. Read its boundary from a point outside all supports. The
chords divide the disk into regions with a tree as \defn{weak dual}, where regions are adjacent
across chords. Root this tree at the region containing that point. For any integer $f\ge2$,
more than $f^f$ edges force either a chain of $f$ edges from the root or a node with $f$
children. In the first case push the corresponding nested chords slightly towards their parent
regions, keeping their ends fixed. In the second case order the child caps along the boundary
and trace arcs in the parent region enclosing respectively the first child, the first two
children, and so on. These arcs have ends at the extreme endpoints of the enclosed intervals.
In either case we obtain $f$ nested caps. Their bounding arcs miss the forest internally,
and each larger cap contains a chord whose interior lies outside the preceding cap.
The arcs may share endpoints.

The boundary interval of a cap avoids the starting point. Give the cap side all its chords
and all pieces whose two society vertices lie in that interval. Put singleton attachments
on the side containing their society vertex. At the two endpoints put the root sets $X,X'$
on both sides, including the entire edgeless attachment at a singleton endpoint.
The piece overlap condition gives separator $X\cup X'$. Nested intervals give nested
separations.

Keep caps with the same ordered pair of endpoint supports and the same separator size $n$.
If the supports differ, each rail on either support has one contact with each separator.
If they coincide, each rail has one or two distinct contacts. The contacts with an inner
separator lie between those with an outer separator along the rail, so their number cannot
increase. Since both separators have size $n$, the count on every rail is the same.
With one contact use the segment between the two cuts. With two, use the two end segments.
This gives $n$ disjoint rail segments joining the first and last separators, allowing
segments of length zero. The number of ordered pairs of supports is bounded, as is $n$ by twice the maximum
depth. Thus we can choose $f$ large enough to keep at least $\thmbound{lem:journal-pump}{\Gamma,k,n}$
caps for some $n$. In the realization given by Lemma~\ref{lem:journal-pump}, only attachment rails are
used in one of the bands. A chord in that band is therefore omitted, contradicting minimality.

For general $S$ induct on $1-\chi(S)$, where $\chi$ denotes the Euler characteristic. A proper
arc between distinct drawing vertices with interior disjoint from the chord forest is
\defn{essential} if cutting along it does not separate off a disk. In particular, a
nonseparating arc is essential. Cut along any such arc. At each endpoint divide its support
and its attachment at the corresponding root set. Duplicate the root set, assigning each
subgraph $K_i$ intact to its side. Retain a singleton attachment on a side with no remaining
subgraph $K_i$. The interval property shows that only the two indicated root sets can acquire
copies, even when both ends are on one support. The vertex increase is at most $2q$ and there
are at most $r+2$ new total supports, still of depth at most $q$. Every boundary circle
still contains a drawing vertex. All drawing edges remain intact. By the splitting observation a
pattern of order at most $k+6q$ suffices. A nonseparating cut decreases $1-\chi$ by one. For a
separating essential cut the two positive values of $1-\chi$ sum to the previous value. Apply
induction to each piece and project. There are at most two pieces, giving a bound in the fixed
parameters.

It remains to consider a surface admitting no such essential arc. It has just one boundary
circle. Indeed the incidence graph of complementary regions and boundary circles is
connected, since regions adjacent across a chord both meet the circle of either endpoint.
If there were two circles, some region would meet two of them. Joining incident drawing
vertices on these circles through the region gives a nonseparating arc.

Every chord is separating, since pushing its interior slightly aside gives an arc
disjoint from the forest with the same ends. Hence no two chords have four alternating
ends on the boundary circle. Place the boundary vertices in the same cyclic order on a
round disk and draw each chord as a straight segment. Distinct segments can meet only
at common ends, so this draws the same forest with the same supports. The disk bound
now applies.
\end{proof}

\begin{lemma}\label{near:vortices}
Fix a compact connected surface $S$ with nonempty boundary, a finite abelian group
$\Gamma$, and nonnegative integers $r,q,k,w$. Let $G$ be a $\Gamma$-labelled graph
presented as a drawing $F$ in $S$ with attachments on at most $r$ supports as above.
Suppose each attachment has depth at most $q$ and treewidth at most $w$. Let $\Pi$
be a pattern in $G$ of order at most $k$, and let $W$ be the union of the societies and
the named vertices of $F$. Deciding realizability reduces to linkage decisions
of order at most $\thmbound{near:vortices}{S,\Gamma,r,q,k}$ in subdrawings of $F-W$ on
the surface obtained by capping the boundary circles of $S$.
\end{lemma}
\begin{proof}
Take $\thmbound{near:vortices}{S,\Gamma,r,q,k}\coloneqq \thmbound{near:chords}{S,\Gamma,r,q,k}$
and enumerate at most this many chords with ends in $W$. For each choose either a
specific drawing edge or its first and last actual edges $xu,vy$, where $x,y\in W$ and
$u,v\notin W$, and the gain required from $u$ to $v$. Distinct chords have disjoint internal
end sets. For $u=v$ require zero inner gain. Test the inside demands in a single drawing query
on $F-W$, reserving every vertex $u=v$. Separately let the attachments, together with the
named drawing vertices as isolated vertices where necessary, receive one edge for each
enumerated chord, with its computed end-to-end gain. Adding the boundedly many chord ends to
all bags in a decomposition of the disjoint attachments gives this graph bounded treewidth.
Test the original pattern there by bounded-treewidth model checking. Any solution in this
second graph expands through the first solution without intersections, since the interior
paths lie outside $W$ and outside the attachments. Conversely, the realization in
Lemma~\ref{near:chords} supplies a successful guess. The enumeration is polynomial since the
bound and group are fixed.
\end{proof}

\section{Reductions in an embedded graph}\label{sec:embedded-reductions}

A \defn{rigid skeleton} is
the simple 3-connected skeleton of a node in the decomposition of a block into
3-connected components. Its edges are called \defn{adjacencies}. An adjacency stands for an edge of
the block or for a subgraph of the block meeting the rest of the skeleton only in the two ends of
the adjacency. A \defn{bundle} at an adjacency consists of parallel edges with the same tail,
the same head and distinct labels. Its members are also called \defn{options}. A \defn{formal edge} is an adjacency designated by the reduction as a placeholder.
It carries a single zero label. At every other adjacency a zero option, like any other
option, records an actual path gain. The
skeleton with its bundles is a labelled multigraph. It need not be a subgraph of the block.

\begin{lemma}\label{thm:linkage-rigid-skeletons}
Let $\Gamma$ be a fixed finite abelian group and let $k$ be a fixed nonnegative integer.
Let $G$ be a $\Gamma$-labelled graph and let $\Pi$ be a pattern in $G$ of order at most $k$.
Deciding realizability reduces to linkage decisions of bounded order on the rigid skeletons
of the blocks of $G$, with every adjacency retained as a nonempty bundle of at most
$|\Gamma|$ options. There is a pair of integers
$(k_0,k_1)\coloneqq \thmbound{thm:linkage-rigid-skeletons}{k}$ such that the block queries
have order at most $k_0$ and the queries on rigid skeletons have order at most $k_1$.
Each resulting skeleton has at most $2k_0+1$ formal edges.
At every other adjacency, each option is the gain of a path in the part of the block
represented by that adjacency. The paths for distinct adjacencies have disjoint interiors
avoiding all other skeleton vertices.
\end{lemma}
\begin{proof}
Work componentwise. We use the same table first on the tree of blocks and
then on the SPQR tree of each block. In either rooted tree, the subgraph
below a node meets the rest only at its \defn{poles}, one or two vertices
belonging to the node. The root has no poles. Mark the poles and every original name in the
subgraph. The table records the linear forests on subsets of the marks
that have a subdivision there avoiding all other marks. Each forest edge
is oriented and carries the gain of its subdividing path. No mark may occur
inside such a path. At the root we test for the prescribed paths, gains and
singletons.

Here is the common step. Call a child \defn{exceptional} if its subgraph
has an original name strictly inside. For pattern order $h$ there are at
most $2h$ such children, since their interiors are disjoint. On the node
itself mark the original names there, its parent poles and the poles of its
exceptional children. Compute the local forest table on these marks.
Take the disjoint union of a local entry and one entry from each exceptional child,
then identify the common marks. Keep the resulting multigraph if it is a linear forest whose ends and
isolated vertices are marked in the whole node network. Suppress the unmarked
degree-two vertices, adding gains along each path. This gives exactly the
table for that network. Indeed, restricting a subdivision gives local entries, and
the shared marks record every possible intersection when gluing them back.

To test a local forest on mark set $M$, guess for each forest edge either a
single actual edge with the prescribed ends and traversal gain or the first and
last actual edges of its path. In the second case require their inner ends to lie outside $M$ and subtract the two
traversal gains from the target. Test the resulting inner demands together,
reserving every mark. Coincident ends of one inner demand require a singleton
of gain zero. Reject shared inner ends of different demands. This tests all
segments on the intact node with at most $2|M|$ demands and reservations.

For the first table root the block--cutvertex tree at a block. Give each other
block the nearest block towards the root as its parent. Each child and its
descendants attach at one cutvertex, which several children may share.
A child with no internal name can supply no edges of a realizing forest,
since a path entering it would have to return through the same cutvertex.
Thus the local query is on the original block, with its inherited drawing
when one is given. Each exceptional child accounts for a name outside the
block, so $|M|\le 2k+1$ and $k_0\coloneqq 4k+2$ suffices.

Handle blocks on at most two vertices directly. For each other block use the
polynomial-time SPQR decomposition of
Gutwenger and Mutzel~\cite{GutwengerMutzelSPQR}. Each node has a skeleton
that is a cycle, a parallel bundle, a single edge or a simple 3-connected
graph. A node of the last kind is \defn{rigid}. A \defn{virtual edge}
represents the network beyond it, which joins its two poles and meets the
skeleton only there. Distinct child networks have disjoint interiors.
A child with no internal name can supply at most one pole-to-pole path,
since two would form a cycle. Replace its virtual edge by the nonempty bundle
of passage gains in its table. Designate each exceptional child edge and the parent
edge, when present, as formal. Forbid them in the end-edge guesses above.
Their ends are marked, so the inner paths already avoid them.

There are at most $2k_0+1$ formal edges. Each exceptional child accounts for
a name outside the skeleton and contributes two marks, so $|M|\le4k_0+2$ and
$k_1\coloneqq 8k_0+4$ suffices. Nonrigid skeletons with their bundles have treewidth
at most two. All tables have bounded size and are computed from the leaves
upwards. The nonformal options at distinct adjacencies come from a real edge
or from paths in distinct child interiors, as required. For an unlabelled
subdivision, a formal edge may likewise use a path in its own network.
\end{proof}

For a graph drawn on a closed surface its \defn{radial graph} is the
bipartite incidence graph of its vertices and complementary regions.
The \defn{radial distance} between two vertices of the radial graph is their distance in that graph.
The vertices of the drawn graph are called \defn{primal} when they are regarded as vertices of
the radial graph.
A \defn{subdrawing} retains its
vertices and edges in their original positions.

For a set $\mathcal F$ of faces of a drawn graph and an integer $r\ge0$, write
$B_r(\mathcal F)$ for the set of drawing vertices within radial distance $r$ of at
least one face in $\mathcal F$.

This section and the next use the following \defn{cutting formula}. If a closed surface
of Euler genus $g$ is cut along $c$ pairwise disjoint two-sided simple closed curves into
$n$ pieces, then
\[
 g=\sum_D g_D+2(c-n+1),
\]
where $D$ ranges over the pieces and $g_D$ is its Euler genus. This follows by
adding Euler characteristics, since the pieces have $2c$ boundary circles in total.
Their adjacency graph is connected, so $c\ge n-1$.

Apply this formula to the boundary curves of disjoint regular neighbourhoods of
connected embedded graphs. Each graph embeds in its capped neighbourhood, so their
minimum Euler genera sum to at most $g$. For a single connected graph of minimum
Euler genus $g$, every complementary piece has genus zero and $c=n-1$. The $n-1$
complementary pieces have $c$ boundary circles between them, one each. Thus every
embedding of the graph in Euler genus $g$ is cellular.

\begin{lemma}\label{cor:mixed-radial-ball-peeling}
Fix a closed surface $\Sigma$, a finite abelian group $\Gamma$, and nonnegative integers
$k,b,r$. Let $G$ be a $\Gamma$-labelled graph drawn in $\Sigma$ and let $\mathcal F$ be
a specified set of at most $b$ faces of this drawing. For a pattern in $G$ of order at
most $k$, deciding realizability reduces to linkage decisions of order at most
$\thmbound{cor:mixed-radial-ball-peeling}{\Sigma,\Gamma,k,b,r}$ on subdrawings avoiding
$B_r(\mathcal F)$.
\end{lemma}
\begin{proof}
We first remove all vertices incident with one specified region,
using its interior point to identify its containing region whenever
edges or vertices are deleted. Induct on the Euler genus $g$ of the surface,
allowing a bounded list of regions at each smaller genus. Work componentwise and then on the
unchanged blocks as in the first part of
Lemma~\ref{thm:linkage-rigid-skeletons}, testing blocks with at most two vertices
directly. In a 2-connected block $B$ let $R$ be the region containing the marked point.
All vertices of $B$ incident with the original region are incident with $R$.
Cut along the boundary of a thin regular neighbourhood of $B$. In the cutting formula
the $n-1$ pieces other than the neighbourhood have $c$ boundary circles in total and at
least one each. Hence the piece inside $R$ has at most $c-(n-2)\le1+\lfloor g/2\rfloor$
boundary circles. They trace the \defn{boundary walks} of $R$ in $B$, together passing
once through each corner of $B$ in $R$.

Suppose some vertex $v$ has two corners in $R$. Join them by an arc through $R$. This gives a
simple closed curve with just $v$ on the graph and edge ends in
both local sectors. A separating such curve would disconnect the
block upon deletion of $v$. The block has no loops and has at least
three vertices, so each sector contains an edge leading to another
vertex. Enumerate visits to $v$ by Lemma~\ref{near:apex}, delete $v$, then cut
along the curve and cap. The surface stays connected and its Euler genus drops by
one for a one-sided curve or two for a two-sided curve. The arc through $R$ cuts
that region into at most two pieces. After deletion and capping their containing
regions still have all surviving vertices incident with $R$ in their closures.
Apply induction to these regions and restore the original drawing in the resulting queries.

Otherwise the boundary walks of $R$ are vertex-disjoint cycles. No vertex
occurs twice, there are no loops, and an immediate reversal of an edge would
require a vertex of degree one. Remove one edge from each cycle to obtain a
path on all its vertices. Use the path as an attachment and leave the omitted
edge in the drawing. The two path edges at each internal vertex border the
same corner of $R$. A thin disk along the path on that side can therefore touch
the remaining drawing at precisely the path vertices, in order. The disks can
be chosen disjoint. This construction also applies to a digon, whose attachment
is a single edge.

Let $S'$ be the surface with the open disks removed, and let $k'$ be the pattern
order in the block. Each attachment has depth and treewidth one. Lemma~\ref{near:vortices}
therefore reduces the problem to linkage queries of order at most
$\thmbound{near:vortices}{S',\Gamma,1+\lfloor g/2\rfloor,1,k'}$ in subdrawings
avoiding all the cycle vertices. These are exactly the vertices of $B$ incident
with $R$. For a bounded list of regions, apply this reduction successively at
their marked points in the original ambient surface.

To remove radial balls, put $X_0\coloneqq \varnothing$ and let $X_j$ contain
the primal vertices at distance at most $2j-1$ from the given faces.
In any subdrawing avoiding $X_j$, every surviving vertex of
$X_{j+1}$ is incident with a region containing a marked centre
point. Indeed the intermediate graph vertices of a shortest radial
path to it have been deleted. The radial incidences can be joined
through their original open regions up to its last vertex. Apply
the preceding reduction to these at most $b$ regions and iterate.
The nesting depth and all orders are bounded at fixed parameters.
\end{proof}

\subsection{Facial cycles in a large embedded graph}

A \defn{rectangular wall} of height $r$ and width $t$ is a subdivision of the following graph on
$[r]\times[2t]$. Take the rectangular grid and delete the edge from $(i,j)$ to $(i+1,j)$
whenever $i\equiv j\pmod 2$. Its horizontal paths are its \defn{rows}. Each pair of columns $2j-1,2j$
contains a vertical path that alternates between them to pass from the first row to the last.
The top vertex in that pair missing its downward edge is its \defn{boundary vertex}. Adding the edge
between the ends of each row before subdivision gives an \defn{annulus wall}. We use its horizontal
cycles and its disjoint vertical paths. A \defn{digon} is a cycle with two edges.

A \defn{disk with ports} is a closed disk bounded by a graph cycle, with specified
distinct vertices on that cycle called its \defn{ports}. An
\defn{enclosure} of such a disk consists of pairwise disjoint graph cycles
$O_1,\ldots,O_m$, disjoint from the boundary cycle of the disk, and
pairwise disjoint paths, one from each port to $O_m$. There are closed
disks $\Delta_1\subset\cdots\subset\Delta_m$ bounded by $O_1,\ldots,O_m$,
with the given disk in the interior of $\Delta_1$ and each $\Delta_i$ in
the interior of $\Delta_{i+1}$. Each path meets the port cycle only at its
start, meets $O_m$ only at its end, and lies in $\Delta_m$ minus the
interior of the given disk. The cycles $O_i$ are called the \defn{surrounding cycles} or \defn{perimeters}, the paths are called
\defn{spokes}, and $\Delta_m$ is the \defn{outer disk} of the enclosure. An
\defn{auxiliary cycle} is a cycle of the graph together with a closed disk that it bounds.
An enclosure of a facial or auxiliary cycle is an enclosure of its disk. The
\defn{facial gain class} of a face with facial gain $a$ is $\{a,-a\}$.

\begin{figure}[htb]
\centering
\setlength{\unitlength}{1pt}
\begin{picture}(282,116)
\linethickness{0.4pt}
\put(18,25){\framebox(150,83){}}
\put(32,35){\framebox(122,63){}}
\put(46,45){\framebox(94,43){}}
\put(69,56){\framebox(48,21){$C$}}
\put(18,66){\line(1,0){51}}\put(117,66){\line(1,0){51}}
\multiput(18,66)(14,0){3}{\circle*{2.3}}
\multiput(140,66)(14,0){3}{\circle*{2.3}}
\put(69,66){\circle*{3.2}}\put(117,66){\circle*{3.2}}
\put(59,71){\makebox(8,8){\small $u$}}\put(119,71){\makebox(8,8){\small $v$}}
\multiput(67,8)(26,0){3}{\line(0,1){37}}
\put(7,8){\line(1,0){170}}
\multiput(67,45)(26,0){3}{\circle*{2.3}}\multiput(67,35)(26,0){3}{\circle*{2.3}}
\multiput(67,25)(26,0){3}{\circle*{2.3}}\multiput(67,8)(26,0){3}{\circle*{2.3}}
\put(180,96){\makebox(96,9)[l]{\small enclosing cycles}}
\put(180,62){\makebox(98,9)[l]{\small paths from the ports}}
\put(180,4){\makebox(97,9)[l]{\small a lower row}}
\end{picture}
\caption{An enclosure of a facial cycle $C$ with ports $u,v$. Three columns run from the
innermost perimeter to a lower row. A dot indicates a contact. In the wall such a contact can
be a shared subpath.}
\label{fig:enclosure}
\end{figure}
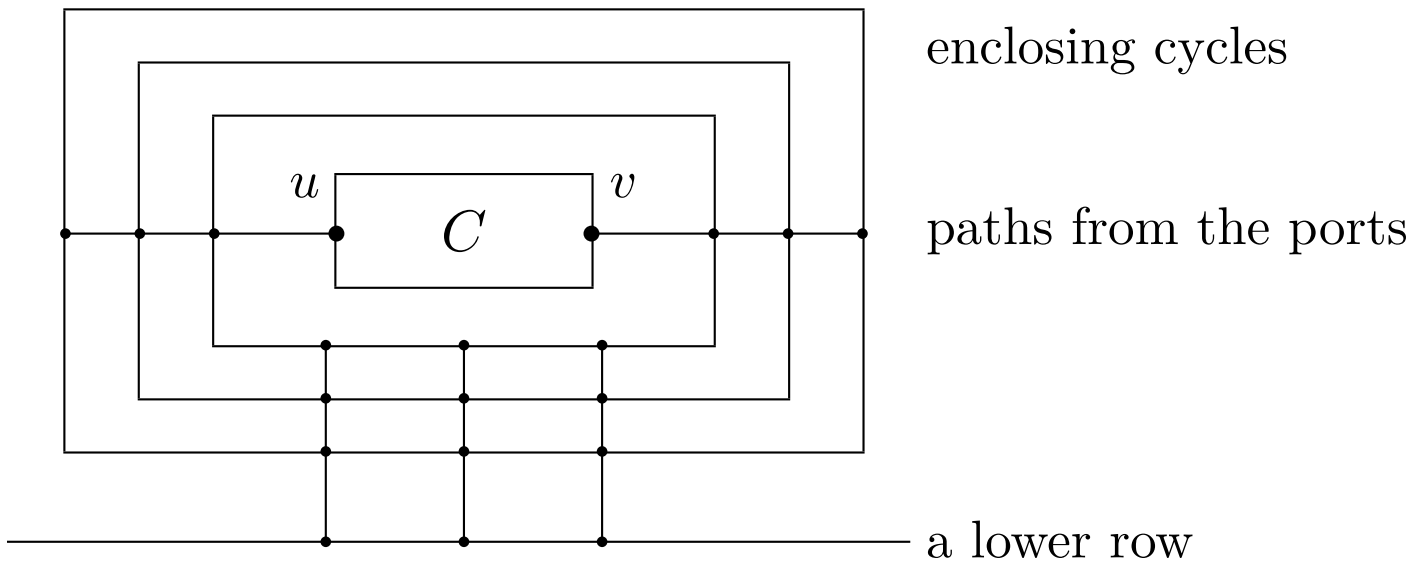

Figure~\ref{fig:enclosure} shows the two kinds of paths used with a facial cycle. The paths
from the ports allow a linkage to reach the cycle. The columns keep its surrounding perimeters
connected to the rest of the wall after a small vertex deletion.

In this section a \defn{vortex disk of adhesion at most $d$} is a closed disk whose
boundary meets the drawing precisely in vertices $x_1,\ldots,x_n$ in cyclic order, with
$n\ge1$ and all $x_i$ distinct.
The graph in the disk must have a path decomposition $(B_1,\ldots,B_n)$ of adhesion at
most $d$ with $x_i\in B_i$ for each $i$. A disk with no graph in it can be omitted.

A \defn{bundle rendition} of a graph $J$ drawn in a surface is a finite family of closed
disks called \defn{cells}. The boundary of each cell meets the drawing in exactly two
vertices, its \defn{poles}. Distinct cells meet only at common poles. Every edge belongs
to a cell and its interior lies in the interior of that cell. The graph in each cell
is the union of one or more paths between its poles, any two meeting only at the poles.
The \defn{torso} has the poles as vertices and one edge per cell joining its poles.
We require the torso to be simple and 3-connected.

For a closed surface $\Sigma$ of positive Euler genus fix a threshold $N_\Sigma\ge3$
from the uniqueness theorem of Seymour and Thomas~\cite{SeymourThomas},
in the form of~\cite[Theorem~5.10.5(a)]{MoharThomassen}. A simple 3-connected graph with
a cellular embedding in $\Sigma$ of representativity greater than $N_\Sigma$ has
minimum Euler genus that of $\Sigma$. The identity on the graph extends to a
homeomorphism between this embedding and every other one in $\Sigma$.
See also Mohar~\cite[Theorems~5.4 and~6.1]{MoharUniqueness}.

\begin{lemma}\label{lem:mixed-represented-digons}
Fix a closed surface $\Sigma$, a finite abelian group $\Gamma$, and positive integers
$m,p,q$. For each integer $s\ge0$ there are nonnegative integer bounds
\[
 (b,d,N)\coloneqq \thmbound{lem:mixed-represented-digons}{\Sigma,\Gamma,m,p,q,s}
\]
with the following property. Let $J$ be a simple 2-connected $\Gamma$-labelled graph
cellularly drawn in $\Sigma$ with simple facial boundaries. Suppose $J$ has a bundle rendition.
On the sphere assume its treewidth exceeds $N$. In positive Euler genus assume
its representativity exceeds $N$.
There are a set $A$ of nonzero facial gain classes and at most $b$ pairwise disjoint vortex
disks of adhesion at most $d$ such that the following hold.
\begin{enumerate}
\item Each face whose gain class is nonzero and outside $A$ lies with its boundary in
one of the vortex disks. More than $s$ vertices of $J$ lie outside each vortex disk.
\item For each class in $A$ there are $q$ two-port enclosures of facial cycles of that
class. There are also $q$ three-port enclosures of auxiliary cycles. The outer disks of all
these enclosures are pairwise disjoint and each enclosure has $m$ surrounding cycles.
For every $X\subseteq V(J)$ with $|X|\le p$, all surrounding cycles disjoint from $X$
lie in one component of $J-X$ with more than $s$ vertices.
\end{enumerate}
The first coordinate $b$ is independent of $s$. On the sphere the second coordinate $d$ is also independent of $s$.
\end{lemma}
\begin{proof}
Give each cell one \defn{colour}, the set of nonzero gain classes of faces incident
with its edges. The number of colours is bounded in terms of $\Gamma$. Each cell joins its
two poles and the torso is simple and 3-connected. Thus the rendition is tight in the
sense of both Bla\v{z}ej, Pilipczuk and Protopapas~\cite[Sections~3--4]{BPPcoarseMenger}
and Paul, Protopapas, Thilikos and Wiederrecht~\cite[Section~4.2]{PPTWindex}.
Their cell-colouring results require the input wall to be \defn{grounded}, meaning
that each of its cycles uses two or more cells. Our walls subdivide simple cubic
3-connected graphs, so every cycle contains at least three vertices of degree three
in the wall. A single cell can contain at most two, its poles.

We use the spherical and positive-genus forms of the cell-colouring theorem
of~\cite[Propositions~4.2 and~10.11]{BPPcoarseMenger}, based
on~\cite[Lemma~8.1]{PPTWindex}. Both forms give a graph assembled from wall segments,
called a \defn{walloid}, together with disks in the given surface. Its base is a cylinder
of height $r$, divided into segments. We use the following features of the output.
\begin{itemize}
\item There is an ordinary $r$ by $r$ wall segment. For each represented colour,
at least $Q$ distinct flap segments each contain a cell of that colour in a closed
disk called a \defn{pocket}. The width of a flap segment is controlled by a parameter $t$.
\item There are at most $b$ vortex segments. Every cell of an unrepresented colour
lies in an inner vortex disk. The corresponding outer disks are pairwise disjoint
and have the ordered path decompositions of adhesion at most $d$ required in (1).
By the disk convention of~\cite[Sections~3 and~4.2]{BPPcoarseMenger}, an outer
boundary meets the drawing only at vertices of the outer cycle of its vortex
segment. In~\cite[Sections~3.3 and~7.1]{PPTWindex} it is the trace of that cycle.
\end{itemize}
The base lies outside each outer vortex disk. Indeed it is connected and misses the
boundary contacts on the outer cycle. The segment's rails end on the base outside
the disk interior by the sandwich condition of~\cite[Sections~3 and~4.2]{BPPcoarseMenger}.

On the sphere the required input annulus wall is a minor of a sufficiently large square
grid (draw concentric cycles with columns between them). Such a grid exists above a
fixed treewidth bound by Robertson and Seymour~\cite{RSV}. Since the wall is
\defn{subcubic}, meaning of maximum degree at most three, a minor model contains a
subdivision. Apply Proposition~4.2 with parameters $r,t,Q$ and the number of colours.
Both $b$ and $d$ are independent of $r$. The output subdivides a 3-connected planar
graph, so Whitney's theorem identifies its drawing with the usual spherical drawing.

In positive Euler genus $g$, put $R\coloneqq a_{10.8}g+2n$ and choose both $r$ and $t$ equal to $R$, with $a_{10.8}$ as
in Lemma~10.8 of~\cite{BPPcoarseMenger}. Take $n\ge4$ large enough for all the
constructions below and take $R>\max\{s+2,N_\Sigma\}$. The input surface
annulus wall for Proposition~10.11 embeds in $\Sigma$ and is subcubic. Thus
Proposition~10.2 of~\cite{BPPcoarseMenger} supplies a subdivision at sufficiently
large representativity. Proposition~10.11 gives the stated output with $b$
depending only on $Q$ and the number of colours, though $d$ can also depend on $R$.

Keep the ordinary, handle and crosscap segments and the horizontal base paths through
the other segments. These paths subdivide the closing edges of an $R$-surface annulus wall.
Increase $n$ to meet the threshold of Proposition~10.1 of~\cite{BPPcoarseMenger}.
By that proposition and Lemma~10.8 there, this subwall has minimum Euler genus $g$
and representativity $R$ in every embedding in $\Sigma$. The whole walloid therefore
has minimum Euler genus $g$ and representativity at least $R$. Every embedding of it in
$\Sigma$ is cellular by the cutting formula.

Draw the base rows concentrically. Use disjoint intervals along the inner rim for
the segment decorations. At the handle and crosscap segments insert the prescribed
handles and crosscaps as in the proof of Lemma~10.8 of~\cite{BPPcoarseMenger}.
Draw the nested flap arcs over their own intervals, and connect each vortex annulus
to its interval in order. The ordinary and flap segments now lie in disks with
disjoint interiors. The walloid subdivides a simple 3-connected graph~\cite[Section~10.4]{BPPcoarseMenger}.
Since $R>N_\Sigma$, the identity on that graph extends to a homeomorphism of the
two embeddings. Restoring the subdivisions carries each segment disk and its flap pocket
to the inherited drawing.

Put $A$ equal to the union of the represented colours. If a face has nonzero
gain class outside $A$, every cell containing a boundary edge lies in an inner
vortex disk. Its connected boundary therefore lies in one such disk. The face is
on the disk side, since the other side contains base vertices on the sphere and
is not an open disk in positive genus. Taking $r>s+2$ puts more than $s$ base
vertices outside each vortex disk. This proves (1).

For (2), take $Q$ larger than $q$ times the number of possible nonzero classes.
For each class in $A$ choose a represented colour containing it and $q$ unused
pockets of that colour. At most $q(|A|-1)$ pockets have been used when the last
class is reached. In each chosen pocket the cell has an edge incident with a
face of the required class. Write $C$ for the boundary of this face.

The flap has two attachment positions for its hosted cell, one per pole~\cite[Section~4.2]{BPPcoarseMenger}.
Its top boundary vertices occur in the order
$v_1,\ldots,v_t$, then $z_1,z_2$, then $u_1,\ldots,u_t$. Above them are nested paths from
$v_j$ to $u_{t-j+1}$, the innermost bounding the pocket with its top
interval~\cite[Section~4.1, Figures~6--8]{BPPcoarseMenger}. Each such path
completes to a cycle by following its end columns down to a common row and
joining them along that row. Follow a column from its top vertex across the
within-column edge if necessary, and stop at its first contact with the closing
row. These instructions also specify the cycles after subdivision.
Moving both columns out and the closing row down gives strictly nested disks in
the local drawing.

Use first the path from $v_{t-1}$ to $u_2$ and the third row. The pocket is
strictly inside this cycle. The chosen face meets the interior near its incident
edge and cannot cross the cycle, so its closure lies in the closed disk.
To obtain $m$ surrounding cycles, put $j\coloneqq t-4(p+1)-10$ with
$t>4(p+1)+2m+30$. For $i=1,\ldots,m$, use the above-row path with left index
$j-i+1$ and close it along row $i+4$. The first disk strictly contains the third-row disk, and each later disk strictly
contains its predecessor. Thus $C$ misses every surrounding cycle. More than $p$ columns run from the bottom of each perimeter
to the lower half of the cylinder.

Only now choose the two ports of $C$. Since $J$ is 2-connected, Menger's theorem
gives two vertex-disjoint paths from $C$ to the outermost surrounding cycle.
Stop each at its first contact with that cycle and keep the portion from its
last contact with $C$. These portions stay in the enclosing disk, and the open
face bounded by $C$ contains no graph. They are the two spokes, with their
distinct ends on $C$ as ports.

For the $q$ auxiliary enclosures take disjoint column intervals in the ordinary
segment. In the upper half of each interval draw one wall rectangle inside $m$
larger rectangles. Make their sides far enough apart to leave more than $p$
disjoint columns from the bottom of each perimeter to the lower half of the
cylinder. At the inner rectangle, taken wide enough, choose three bottom vertices at downward edges
and follow their columns to the outermost perimeter. The outer disks are strictly
inside their respective segments or intervals, hence pairwise disjoint.

Finally take more than $p$ lower rows, each of length greater than $s$. For any
deletion of at most $p$ vertices, some such row survives intact. Every surviving
perimeter reaches it along one of its more than $p$ disjoint columns, proving (2).
On the sphere only $r$ needs to grow with $s$. In positive genus the common
dimension $R$ can meet all the height and width requirements simultaneously.
\end{proof}

For nonnegative integers $s,d$, put
$\thmbound{lem:confined-face-regions}{s,d}\coloneqq(2d+1,s+2d+1)$.

\begin{lemma}\label{lem:confined-face-regions}
Let $\Sigma$ be a closed surface, let $\Gamma$ be a finite abelian group, and let $J$ be a
connected $\Gamma$-labelled graph cellularly drawn in $\Sigma$.
Let $b,d,s$ be nonnegative integers. Let $A$ be a set of facial gain classes and
$\mathcal D$ a family of at most $b$ pairwise disjoint vortex disks of adhesion at most $d$.
Suppose every face of nonzero gain class outside $A$ lies with its boundary in a disk
of $\mathcal D$, and more than $s$ vertices of $J$ lie outside each disk. Let $K\le\Gamma$
be generated by the elements of the classes in $A$ and put
$(r,v)\coloneqq\thmbound{lem:confined-face-regions}{s,d}$. There is a set $\mathcal F$
of at most $b$ faces such that each component of $J-B_r(\mathcal F)$ either lies in a disk
of $\mathcal D$ or has a labelling modulo $K$ that extends flatly off at most $b$ points of
$\Sigma$. All components in these disks are planar. If $J$ is
$(s,c)$-unbreakable for some integer $c\ge r$, each component in a disk has treewidth
at most $v$.
\end{lemma}
\begin{proof}
Consider a vortex disk $D$ with at least two boundary contacts $x_1,\ldots,x_n$ and
ordered bags $B_1,\ldots,B_n$. In the gap from $x_n$ to $x_1$ choose a point $a$. As centre
use the original face containing it. For $1\le j<n$ choose a point $a_j$ in the gap following
$x_j$. We show that the face containing $a_j$ is at radial distance at most $2d$ from the
centre. Its contacts $x_j,x_{j+1}$ are then at distance at most $2d+1$ from the centre.

Deleting $S\coloneqq B_j\cap B_{j+1}$ separates the contacts in the two open boundary arcs
from $a$ to $a_j$ in $J\cap D$, since the bags containing a vertex are consecutive. Hence
there is an arc in $D$ from $a$ to $a_j$ missing the surviving graph. To see this, thicken
one of the two boundary arcs and all surviving components meeting its relative interior. They
miss the opposite open arc. The frontier of a sufficiently small regular neighbourhood on the
interior side of the boundary contains the required arc after shortening its ends. Perturb this
arc to avoid vertices and to cross edges transversely. Every crossed edge has an end in
$S$ incident with both faces on its sides. Reading the crossings gives a radial walk
between the two gap faces whose primal vertices all lie in $S$. Removing repetitions
gives a path of length at most $2|S|\le 2d$. A single
contact needs only an incident face, and an empty boundary needs no centre.

Write $\mathcal F$ for the chosen faces and $X\coloneqq B_{2d+1}(\mathcal F)$. This set contains every boundary contact. A surviving component
meeting an open disk therefore lies entirely in that disk. Let
$S_i\coloneqq \{x_i\}\cup(B_{i-1}\cap B_i)\cup(B_i\cap B_{i+1})$,
taking empty end intersections. Each vertex of $B_i-S_i$ occurs
only in this bag and is strictly inside the disk. Its neighbours
all lie in $B_i$. Since the graph has more than $s$ vertices
outside the disk, unbreakability and $|S_i|\le 2d+1\le c$ bound
$|B_i-S_i|$ by $s$. This gives the asserted width bound, inherited
by subgraphs.

Consider an exterior component $C$ of $J-X$. Each deleted primal ball is connected, since
consecutive facial boundary walks along a shortest radial path share a vertex.
Every complementary region of $C$ that is not an old face contains a vertex of $X$.
Indeed such a region contains an omitted vertex or edge of $J$. Since $C$ is induced,
an omitted edge has an end outside $C$. By connectedness of $J$, each component
of $J-V(C)$ meets $X$, as an edge from it to $C$ must have its outside end in $X$.
Thus each such region contains a whole deleted ball, and there are at most $b$ of them.

Puncture each of these regions once. Every other region is an old face with boundary
in $C$. This boundary misses all vortex disks, so the face has gain zero modulo $K$
by hypothesis. The proof of Lemma~\ref{thm:find-flatness-punctures} extends the
labelling region by region, with one puncture sufficient in each exceptional region. Both the
width bound and the extension persist on subdrawings.
\end{proof}

\section{Routing in an embedded graph}\label{sec:surface-routing}

\subsection{Facial cycles of a rigid skeleton}
Throughout this section the ambient surface is closed unless holes are specified. The
existence of a noncontractible simple closed curve with at most $n$ intersections implies that deleting at
most $n$ vertices leaves each component embeddable in smaller Euler genus. Include one
endpoint of each crossed edge, delete the vertices on the curve, cut along the curve, and cap
the new boundary components.

For a skeleton with labelled bundles, a \defn{bundle encoding} is the incidence encoding
of the labelled multigraph obtained by replacing each adjacency by its bundle. A
\defn{reference option} is one chosen member of a bundle. For a subgroup $L\le\Gamma$,
two options are \defn{congruent modulo $L$} if their gains in the common direction have
the same image in $\Gamma/L$, and \defn{incongruent} otherwise.

\begin{proposition}\label{prop:surface-spqr-subgroup-guard}
Fix a closed surface $\Sigma$, a finite abelian group $\Gamma$, a subgroup $L\le\Gamma$,
and nonnegative integers $k,h$. Let $G$ be a $\Gamma$-labelled graph drawn in $\Sigma$
whose labelling modulo $L$ has a flat extension off at most $h$ punctures. Let $S$ be
a rigid skeleton with the bundles produced by Lemma~\ref{thm:linkage-rigid-skeletons} for a
pattern in $G$ of order at most $k$. There are nonnegative integer bounds
$N\coloneqq \thmbound{prop:surface-spqr-subgroup-guard}{\Sigma}$ and
$H_0\coloneqq \thmbound{prop:surface-spqr-subgroup-guard}{\Gamma,k,h}$, each written with the
arguments on which it depends, for which at least one
of the following holds.
\begin{enumerate}
\item Deleting at most $N$ vertices of $S$ leaves only components of Euler genus
strictly less than that of $\Sigma$.
\item The bundle encoding of $S$, with one reference option marked at each adjacency,
satisfies a monadic second-order condition depending only on $\Sigma,\Gamma,L,k,h$.
From any structure satisfying this
condition one can construct in polynomial time a cellular drawing of its incidence graph
in $\Sigma$ with at most $H_0$ defective faces modulo $L$. The incidence graph is simple
and 2-connected, and the drawing has simple facial boundaries.
\end{enumerate}
On the sphere the second conclusion always holds.
\end{proposition}
\begin{proof}
Write $g$ for the Euler genus of $\Sigma$. On the sphere put $N\coloneqq 0$ and
use Whitney's theorem for the 3-connected skeleton. For $g>0$ put
$N\coloneqq N_\Sigma$ with the uniqueness threshold of Section~\ref{sec:embedded-reductions}.
Test all deletions of at most $N$ vertices of $S$.
A successful deletion also works for its incidence graph. In each remaining component
draw the parallel options along their reference edge, subdivide them, and restore the leaves
at adjacencies with one deleted end.

If the test fails, even the empty deletion fails. Thus $S$ has minimum Euler genus $g$,
and every embedding in $\Sigma$ is cellular. Each has representativity greater than $N$
by the cutting argument above. Uniqueness therefore fixes the facial walks on the
given graph.
In either genus
the faces are bounded by simple cycles, and distinct facial boundaries meet in at most a
vertex or an edge~\cite[Proposition~5.5.12]{MoharThomassen}.

With $(k_0,k_1)\coloneqq \thmbound{thm:linkage-rigid-skeletons}{k}$, there are at most
$a\coloneqq 2k_0+1$ formal edges in one skeleton. Each simultaneous choice of one option per
adjacency can be drawn as a subdivision of $S$ in the original block. The chosen pole-to-pole
paths occupy the distinct sides of the skeleton edges, which meet only at poles. For a formal
edge choose any path in its side, including the parent side when required. Every face of this
subdivision avoiding the punctures has quotient gain zero. Assigning the formal gains to the
formal edges alters at most two facial gains per formal edge. Consequently every simultaneous choice
has at most $D\coloneqq h+2a$ nonzero facial gains. If an adjacency admits two incongruent options,
choose between two such options equiprobably and independently at all such adjacencies.
Conditional on the other choices, each incident face has probability at least one half of a
nonzero gain, since its boundary uses the adjacency once. Hence at most $2D$ faces are
incident with an incongruent adjacency. Each edge has two distinct incident faces and each
pair of faces shares at most one edge. At most $b\coloneqq \binom{2D}{2}$ adjacencies admit incongruent
options.

In the incidence encoding, mark one reference edge vertex per pair of adjacent original
vertices. Compare all gains at that pair in the direction of its reference, negating a stored
label if necessary. The condition requires this marked bundle encoding, a simple
3-connected reference graph embeddable in $\Sigma$, the failed deletion test when $g>0$,
at most $D$ nonzero reference faces, and at most $b$ incongruent adjacencies. Each bundle has
at most $|\Gamma|$ distinct options. To define facial cycles, subdivide every edge of a
proposed cycle once and join the subdivision vertices to a new vertex. The enlarged graph
embeds in $\Sigma$ exactly when the cycle is facial. In the reverse direction all open spokes
lie in one face of the reference embedding. Its simple boundary contains the entire cycle.
Fixed-surface embeddability is monadic second-order definable by the finite obstruction
theorem of Robertson and Seymour for minor-closed graph classes~\cite{RSXX}. The test with the added vertex can be
interpreted in the incidence structure, using the existing edge vertices as subdivision
vertices. A bounded deletion and its components can be quantified, as can the facial gain
using the fixed finite group. This proves definability. To construct the drawing, compute an
embedding of the reference graph on the indicated surface, using the fixed-surface embedding
algorithm of Mohar~\cite{MoharEmbedding}. Draw the other options parallel in edge strips and subdivide
each option. Store its gain on the first segment from its tail, with zero on the other.
Each original face follows the outermost option on its side of each strip, so its boundary
is still a simple cycle after subdivision. The other faces are subdivided digons between
consecutive options.
Congruent strips introduce only zero quotient digons. The incongruent strips account for at
most $(|\Gamma|+2)b$ additional defective faces, including their exterior neighbours.
Thus we may take $H_0\coloneqq D+(|\Gamma|+2)b$. Puncturing each defective face gives the required extension.
\end{proof}

\subsection{Prescribing routes between faces}
For tangles we assign the edges of a separation to its two subgraphs, so their union is the
graph and their edge sets are disjoint. A \defn{tangle of order $\theta$} chooses one \defn{small side}
of each separation of order below $\theta$, subject to two conditions. No three chosen
small subgraphs cover the graph, and no chosen small side contains all its vertices. For a
tangle of order $\theta$ in a cellular drawing, a vertex set $Z$ of size at most $\theta$ is
\defn{free} if no small side of order less than $|Z|$ contains $Z$. A tangle is
\defn{respectful} if every noose meeting fewer than $\theta$ vertices bounds a disk whose side
of the induced separation is small. The radial graph has the drawing vertices and one point
per face as its vertices. For each incidence of a vertex with a face there is an edge through
that face to its point, with one edge per corner. The \defn{inside} of a radial cycle of length less
than $2\theta$ is the closed disk selected by respectfulness. The inside of a closed radial
walk of such length consists of its drawing and the insides of all simple cycles in its
support. The \defn{tangle distance} of two distinct faces is the least half-length of a closed radial
walk of length less than $2\theta$ whose inside contains both face points, or $\theta$ if
there is no such walk. A \defn{hole} is a closed disk meeting the drawing only at a nonempty
set of vertices on its boundary. These vertices are its \defn{sites}. Its open disk lies in a
face, which is used for measuring distances. For a family $\mathcal H$ of disjoint holes,
write $Z(\mathcal H)$ for the union of its site sets. A partition $\mathcal P$ of
$Z(\mathcal H)$ is \defn{topologically feasible} if there are pairwise disjoint closed
disks $D_P$ in the surface, indexed by $P\in\mathcal P$, with
$D_P\cap\bigcup\mathcal H=P$ for each part $P$.

We use the following result of Robertson and Seymour~\cite[(3.2)]{RSXII}, in the formulation
of Liu and Thomas~\cite[Theorem~5.9]{LiuThomas}, which applies in particular to the sphere. Taking the maximum of its thresholds
over smaller numbers of disks and sites gives the nondecreasing bound below.
\begin{theorem}\label{suff:RS}
Let $\Sigma$ be a closed surface and let $h,z$ be nonnegative integers. There is a nondecreasing (in $h$ and $z$) function $\thmbound{suff:RS}{\Sigma,h,z}$ with the following property. Let $G$ be a graph cellularly drawn in $\Sigma$
with a respectful tangle of order at least $\thmbound{suff:RS}{\Sigma,h,z}$. Let $\mathcal H$ be a family of at most
$h$ pairwise disjoint holes with at most $z$ sites in total. Suppose the sites of each
hole are free and the holes lie in distinct faces at pairwise tangle distance at least $\thmbound{suff:RS}{\Sigma,h,z}$.
Let $\mathcal P$ be a topologically feasible partition of $Z(\mathcal H)$. Then $G$
contains pairwise vertex-disjoint connected subgraphs $G_P$, indexed by $P\in\mathcal P$,
such that $V(G_P)\cap Z(\mathcal H)=P$ for each $P\in\mathcal P$.
\end{theorem}

For a finite set $Z$ of marked points in a surface $M$ and a positive integer $e$, write
$H_1(M,Z)$ for \defn{first relative homology} with coefficients in $\mathbb Z/e\mathbb Z$.
Concretely, triangulate $M$ with $Z$ among its vertices. Take oriented edge chains whose
boundary is supported on $Z$ and identify chains differing by sums of oriented triangle
boundaries. Subdivision gives the same group. An oriented arc with ends in $Z$ defines such a
class. For paths with the same ends, their difference is an absolute cycle. A \defn{star} at a disk
consists of a vertex in the interior and disjoint arcs from it to the marked boundary points,
meeting only at that vertex. 

Consider a closed surface $\Sigma$ with disjoint closed disks of two kinds, which
together we call the \defn{marked disks}. A \defn{root
disk} has a nonempty finite set of individually marked boundary sites. A \defn{two-port
disk} has two individually marked boundary points, its ports. An \defn{oriented pairing}
is a set of ordered pairs of distinct root sites in which no site occurs twice. Each
pair may be assigned a list of two-port disks, with no disk repeated in or between lists.
Write $Z$ for the root sites and $M$ for the surface obtained from $\Sigma$ by deleting
the open root disks. For an integer $e\ge1$, a tuple of classes in $H_1(M,Z)$ modulo $e$,
one per pair, is \defn{attainable} if disjoint proper arcs in $M$ join the pairs in their
specified directions and have those classes.

A \defn{frame} is a finite graph embedded in $\Sigma$ meeting each marked disk exactly in
a star with a spoke to every site or port. Its \defn{exterior graph} is obtained by removing
the interiors of the stars and retaining their sites and ports. It is loopless. Each site
or port has one exterior edge leading to a vertex outside the marked disks, and every
vertex outside the disks has degree at most three. Another embedding of this abstract
graph \defn{respects the marked disks} if its intersection with each closed disk is the
prescribed star and it sends every site and port to its individually specified point.

\begin{lemma}\label{prop:surface-relative-frame}
Let $\Sigma$ be a closed surface and let $e$ be a positive integer. Let $n_1,n_2,n_3$
be nonnegative integers. Give $\Sigma$ at most $n_1$ root disks with at most $n_2$ sites
in total, at most $n_3$ two-port disks, an oriented pairing and assigned lists.
There is a frame $F$ for these data with
\[
 |V(F)|+|E(F)|\le \thmbound{prop:surface-relative-frame}{\Sigma,e,n_1,n_2,n_3}
\]
and with
the following property. In every embedding of $F$ respecting the marked disks, every
attainable tuple in relative homology modulo $e$ is realized by vertex-disjoint paths
of $F$ joining the prescribed pairs in their given directions. These paths avoid the open
root disks and unused root sites, visit their assigned two-port stars in order, and each avoids
every other two-port disk.
\end{lemma}
\begin{proof}
There are finitely many tuples, since $M$ has a finite triangulation. Choose one family of
piecewise linear arcs per attainable tuple. Each family can be made to visit its assigned
two-port disks in order. Place tiny disks astride the arcs at the desired visits and move them
to the given disks by an ambient isotopy fixing the root boundaries. To do this shrink and
move their centres in succession avoiding the other centres, then expand and adjust the two
port positions in the collars. Both bijections on two points of a circle are realized by
orientation-preserving homeomorphisms. Put unused disks off the arcs. Replace each crossing of
a disk by its two spokes. These modifications preserve the classes measured after filling all
the two-port disks.

Arrange common short exterior stubs at each port or site across all representatives, including
an unused stub where needed. Take a slightly larger disk, called a \defn{cap}, around each marked
disk so the part in its annular collar consists of the stubs. Outside the caps put the
finitely many representatives in general position except on the common stubs. Subdivide
intersections and extend the resulting graph and the cap rims to a triangulation of the
complement. Subdivide outside the caps so each ordinary triangle meets a cap in at most one
edge or vertex. Regard the caps as additional polygonal faces. Truncate each vertex to a small
cycle in the cyclic order of its incident edges. The resulting cubic graph $B$ has simple
facial cycles and distinct facial closures meet in at most one edge or vertex. Thus the
simple graph $B$ is 3-connected~\cite[Proposition~5.5.12]{MoharThomassen}.

On the sphere Whitney's theorem fixes the facial cycles of $B$. In positive genus
let $N\coloneqq\thmbound{prop:surface-spqr-subgroup-guard}{\Sigma}$. We can ensure
representativity greater than $N$ by refining before truncation. Collapse each cap to
a point. The quotient is homeomorphic to $\Sigma$, and the collapse is a homotopy
equivalence, as one sees in disjoint collars. Choose a metric on the quotient and
a finite cover by coordinate disks with Lebesgue number $\delta$. Refine the triangulation
and choose the truncation cycles small enough that the image of each closed face has
diameter less than $\delta/(N+1)$. A simple closed curve meeting $B$ in at most $N$
points lies in a chain of at most $N$ such faces or in a single face. Its image lies
in a coordinate disk, so the curve is contractible. The choice of $N$ now gives minimum Euler genus
and uniqueness as in Proposition~\ref{prop:surface-spqr-subgroup-guard}.
In either genus the identity on $B$ extends over the facial disks to a homeomorphism
between its given drawing and any other embedding in $\Sigma$.

At each original end of a stub on a cap rim, the truncation creates a rim edge between the two
old rim-edge incidences. Subdivide this new edge and attach the stub there. A representative
can continue along the vertex cycle to its original next incidence. It uses each ordinary old
vertex at most once within its disjoint family, so the local replacements preserve
disjointness and classes. For a cap with only one site, subdivide a spare rim edge, connect
the two subdivisions by a chord, and attach the one stub at the midpoint of that chord. Its
representatives follow only the half-chord to the original attachment. All exterior degrees
are at most three.

In every embedding, uniqueness identifies the subdivided $B$ by a homeomorphism. A component
attached at interior points of two distinct edges of one cap rim lies in that cap, since no
other face has both edges. This applies to each star of size at least two and to each added
chord. The corresponding marked disk misses $B$ and meets the star, so is inside the same
cap. Inside the annular collar the disjoint arcs connecting corresponding marked points are
unique up to homeomorphism fixing the outer rim and respecting the individual sites. Cut along
the arcs and extend boundary maps over disks to see this, ignoring the unused half-chord in a
one-site cap. Thus the new embedding of the graph used by the representatives is a
homeomorphic image respecting the marked disks and points. The induced map $h_*$ of such
a homeomorphism $h$ permutes the jointly attainable tuples. For a requested tuple $x$, the images of the paths for
$h_*^{-1}(x)$ realize $x$ in the new embedding. Up to homeomorphism there are finitely many choices of the
marked disks, the cyclic orders of their sites, the pairing and its lists within the
given bounds. For each type fix a frame constructed above. The maximum of their sizes
gives the stated bound and is nondecreasing in $n_1,n_2,n_3$.
\end{proof}

A \defn{correction digon} is a facial cycle of length two whose two vertices are
specified as ports. A list of correction digons \defn{covers} a subgroup $L$ if its facial
gains belong to $L$ and their subset sums fill $L$. Either orientation of each boundary
cycle may be used, since changing a sign translates the set of subset sums by a member of
$L$. For the trivial subgroup the empty list suffices. If a set generates $L$, then some list
of at most $|L|-1$ of its elements, with repetitions allowed, has subset sums filling $L$.
Start with the subset sums $S\coloneqq \{0\}$. If $S\ne L$, some generator $a$ has $S+a\ne S$,
since invariance under every generator would make $S$ invariant under $L$. Appending $a$
replaces $S$ by $S\cup(S+a)$, increasing its size. At most $|L|-1$ steps suffice.

Let $G$ be a $\Gamma$-labelled graph cellularly drawn in a closed surface $\Sigma$ and let
$L\le\Gamma$. A \defn{routing system} for $L$ in $G$ consists of the following data.
\begin{enumerate}
\item Pairwise disjoint holes, called \defn{root holes}, and a pattern $\Pi$ in $G$ all of
whose named vertices are sites of root holes.
\item A flat extension of the labelling modulo $L$ off punctures in the interiors of the root
holes. We write $M$ for $\Sigma$ minus the open root holes and use this extension
to evaluate the quotient gain of arcs in $M$.
\item For each ordinary demand of $\Pi$, an ordered list of correction digons covering $L$.
No digon occurs twice within or between these lists.
\item A set of facial triangles, called \defn{auxiliary triangles}. They provide three sites
each for representing the other vertices of a frame.
\end{enumerate}
The closed faces bounded by the correction digons and by the auxiliary triangles are
required to be pairwise disjoint and to miss the root holes. In each of these faces we fix a
hole whose sites are all the vertices of its boundary. These holes and the root holes are
the \defn{holes of the system}. The system has \defn{size at most $(n_1,n_2,n_3)$} if it has
at most $n_1$ root holes, at most $n_2$ root sites in total, and at most $n_3$ correction
digons in total. For an integer $T$, the system is \defn{$T$-spread} if $G$ has a respectful
tangle of order at least $T$ in which the site set of each hole of the system is free, and
the holes of the system lie in distinct faces at pairwise tangle distance at least $T$.

\begin{lemma}\label{cor:free-surface-frame}
Let $\Sigma$ be a closed surface, let $e$ be a positive integer, and let $n_1,n_2,n_3$ be
nonnegative integers. There is a pair of integers
\[
 (j,T)\coloneqq \thmbound{cor:free-surface-frame}{\Sigma,e,n_1,n_2,n_3}
\]
with the following property. Let $\Gamma$ be a finite abelian group
and let $L\le\Gamma$ with $\exp(\Gamma/L)=e$. Let $G$ be a $\Gamma$-labelled graph cellularly
drawn in $\Sigma$ with a $T$-spread routing system for $L$ of size at most $(n_1,n_2,n_3)$
with $j$ auxiliary triangles, pattern $\Pi$ and complement $M$ of the open root holes.
\begin{enumerate}
\item The graph $G$ contains a subdivision of the exterior graph of a frame supplied by
Lemma~\ref{prop:surface-relative-frame} for the following data. The root disks are the root
holes, the two-port disks are the closed faces of the correction digons, the oriented pairing
consists of the ordinary demands of $\Pi$, and the lists are those of the system. The
subdivision respects the root sites and the correction ports and uses no edge of a
correction digon.
\item If disjoint proper arcs in $M$ join the ordinary demands of $\Pi$ in their directions
with their prescribed quotient gains, then $\Pi$ has a realization in $G$. Each of its paths
uses one edge of every digon in its list, in the order of the list.
\end{enumerate}
\end{lemma}
\begin{proof}
Put
\[
 j\coloneqq \thmbound{prop:surface-relative-frame}{\Sigma,e,n_1,n_2,n_3},\qquad
 T\coloneqq \max\{n_2+4,\thmbound{suff:RS}{\Sigma,n_1+n_3+j,n_2+2n_3+3j}\}.
\]
Both bounds are nondecreasing in the three counts. Choose the frame for the given marked
disks, pairing and lists from Lemma~\ref{prop:surface-relative-frame}. Its exterior graph
has at most $j$ vertices other than marked sites or ports. To represent each such vertex,
replace it by a little disk with three ports, one for each incident edge and the rest unused.
Move these disks to distinct auxiliary triangles by a homeomorphism fixing the marked disks.

A thin disk along each exterior edge pairs its end sites. Give each unused site a singleton
part. These disks and small disks at the singletons witness topological feasibility.
Theorem~\ref{suff:RS} gives pairwise disjoint paths joining the pairs and avoiding every
other site. Each pair meets each vertex triangle in at most one site because the exterior
graph is loopless. It meets each correction digon in at most one site because a correction
port's exterior edge leads to an ordinary vertex. Thus no connecting path uses a boundary
edge of either kind of face. On each triangle take one used port as the vertex image and join
it by a triangle edge to every other used port. For an isolated vertex take any port.
This completes the subdivision in (1).

For (2), take any rooted subdivision of the chosen frame's exterior as in (1) and insert
the prescribed stars. Lemma~\ref{prop:surface-relative-frame} gives disjoint routes
with the relative classes of the specified arcs. These routes have the required quotient gains.
Indeed, two paths with the same ends and relative class differ by a cycle that vanishes in
$H_1(M;\mathbb Z/e)$, since $H_1(Z)=0$ for the finite set $Z$ of root sites. Their gain
difference vanishes because the quotient transport factors through this homology group.
Replace each used two-port star by either edge of its digon. The replacement is inside an
unpunctured disk, so preserves the quotient gain. Switching that edge changes the full gain
by the facial gain up to sign. Each list covers $L$ even after these signs are chosen.
The lists use disjoint disks, so the switches can realize the required gains independently.

Figure~\ref{fig:gain-switches} illustrates the choices for a list of two digons.

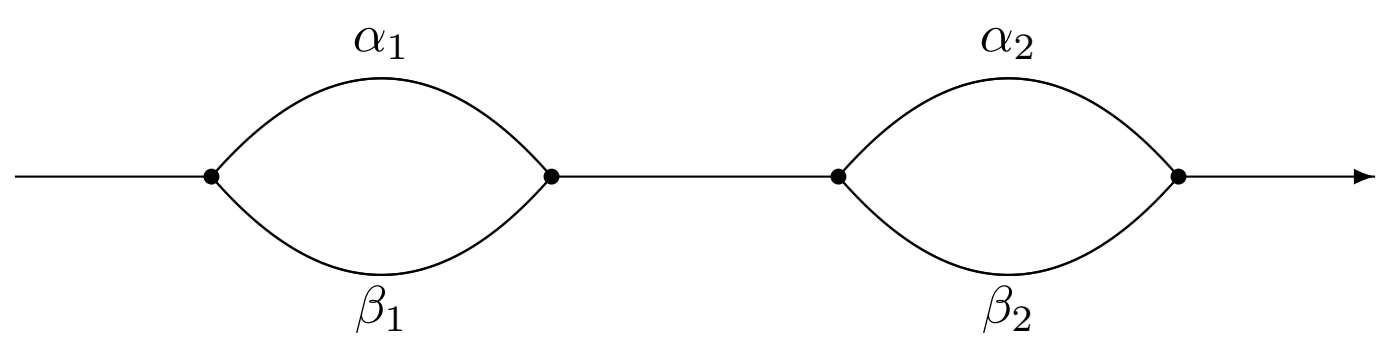
\begin{figure}[htb]
\centering
\setlength{\unitlength}{1pt}
\begin{picture}(280,67)
\thinlines
\put(12,31){\line(1,0){37}}
\put(113,31){\line(1,0){54}}
\put(231,31){\vector(1,0){37}}
\qbezier(49,31)(81,68)(113,31)
\qbezier(49,31)(81,-6)(113,31)
\qbezier(167,31)(199,68)(231,31)
\qbezier(167,31)(199,-6)(231,31)
\put(49,31){\circle*{3}}\put(113,31){\circle*{3}}
\put(167,31){\circle*{3}}\put(231,31){\circle*{3}}
\put(81,51){\makebox(0,10){\small $\alpha_1$}}
\put(81,1){\makebox(0,10){\small $\beta_1$}}
\put(199,51){\makebox(0,10){\small $\alpha_2$}}
\put(199,1){\makebox(0,10){\small $\beta_2$}}
\end{picture}
\caption{Two correction digons on one path. Read each arc from left to right.
Replacing the lower arc of digon $i$ by its upper arc adds the facial gain
$a_i\coloneqq \alpha_i-\beta_i$. If the path using both lower arcs has gain $b$, the
four choices give $b$, $b+a_1$, $b+a_2$ and $b+a_1+a_2$.}
\label{fig:gain-switches}
\end{figure}

\end{proof}

\subsection{Separating the terminal faces}
For a pattern in a drawn graph, a \defn{site incidence} on a face is a specified
boundary vertex of that face. We count each vertex--face pair once. A system of
\defn{designated faces} modulo a subgroup $L$ includes every defective face modulo $L$.
Each named vertex is assigned a site incidence on a designated face, and each designated
face is given at least one site incidence. Extra sites need not be ends of demands.

For a graph $J$, an integer $s\ge0$ and $X\subseteq V(J)$, write $C_X$ for the
unique component of $J-X$ with more than $s$ vertices whenever it exists. Let $H$ be a
subgraph of $J$, let $U\subseteq V(J)$ and let $\theta$ be a positive integer.
Suppose $C_{U\cup S}$ exists and lies in $H-S$ for every separator $S$ of $H$ of order
less than $\theta$. Write $\mathcal T_H$ for the collection that, at each such separation,
selects the side disjoint from $C_{U\cup S}$.

\begin{lemma}\label{lem:genus-carrier-normalization}
Fix a closed surface $\Sigma$ of Euler genus $g$, a finite abelian group $\Gamma$, a
subgroup $L\le\Gamma$, and nonnegative integers $z,w_0$. Let $T$ be a nondecreasing
integer-valued function on the nonnegative integers with $T(w)\ge\max\{3,w\}$ for all $w$.
There are positive integers
\[
 (w_*,u,\theta,c)\coloneqq \thmbound{lem:genus-carrier-normalization}{T,z,w_0}
\]
with the following property. Let $s$ be a nonnegative integer, let $J$ be a connected
$\Gamma$-labelled graph drawn in $\Sigma$, and let $\Pi$ be a pattern in $J$ with at most
$z$ designated faces and at most $w_0$ site incidences as above, modulo $L$. Suppose
$J$ is $(s,c)$-unbreakable and has more than $3s+c$ vertices. If $g>0$, suppose also
that no deletion of at most $s+c$ vertices leaves only components of Euler genus less
than $g$.

For every realization of $\Pi$ there are a connected induced subgraph $H$ of $J$,
a set $U\subseteq V(J)$ with $|U|\le u$, and a residual pattern in $H$ with designated
faces and site incidences satisfying the following conclusions.
\begin{enumerate}
\item The given realization splits into fixed edges, each with an end outside $H$, and a
realization of the residual pattern. Every realization of that pattern pastes to the recorded fixed
edges to realize $\Pi$ in their prescribed order.
\item The inherited drawing of $H$ is cellular. Its labelling modulo $L$ extends flatly
off one puncture in each designated face. All ends and reservations of the residual
pattern have designated site incidences. If $w$ is their total number of site incidences,
including additional sites of designated faces, then $w\le w_*$. The residual pattern
has order at most $w_*$. The collection $\mathcal T_H$ exists and is a respectful tangle of
$H$ of order $\theta$. In $\mathcal T_H$ the site set of each designated face is free and
the designated faces have pairwise distance at least $T(w)$.
\item The graph $H$ omits at most $s+u$ vertices of $J$. The graph $H-U$ is a union of
components of $J-U$ with all their edges, including the unique component of $J-U$ with
more than $s$ vertices.
\end{enumerate}
\end{lemma}
\begin{proof}
For each $X$ with $|X|\le c$ let $C_X$ be the unique component of $J-X$ of more than $s$
vertices. For existence, if all components had size at most $s$, accumulate them until their
union first exceeds $s$. Its size is at most $2s$ and the remaining components have more than
$s$ vertices, violating unbreakability. The same condition gives uniqueness and
$C_{X'}\subseteq C_X$ for $X\subseteq X'$. For $g>0$, some component of $J-X$ still has Euler
genus $g$. There is at most one, because disjoint connected embedded graphs have sum of
minimum Euler genera at most $g$ by the cutting formula of Section~\ref{sec:embedded-reductions}. Its size exceeds $s$, since deleting it as well would give the prohibited
smaller-genus deletion. Thus it is $C_X$.

Start with $U\coloneqq\varnothing$ and $H\coloneqq J$. Throughout, $H$ is a connected induced
subgraph of $J$ with the component property of~(3). With
$|U|+3\theta\le c$, orient every separation of $H$ of order less than $\theta$ towards the
component containing $C_{U\cup S}$, where $S$ is its separator. The common component for the
union of three such separators avoids all three small sides. This proves the tangle axioms. The
connected graph $H$ is cellular on the sphere. At positive genus it
contains $C_U$, so it too has minimum Euler genus $g$ and is cellular. Any short noose misses $C_{U\cup S}$. On the sphere select its disk away from
that component. In positive genus the component's minimum Euler genus rules out a
noncontractible noose by cutting and capping, and places the component outside the unique disk
of a contractible noose. This proves respectfulness.

If two designated faces have distance below $T(w)$, take a witnessing radial walk and let $P$
be its connected image and $A$ its set of primal vertices. Keep the component $H'$ of $H-A$
containing $C_{U\cup A}$, and add $A$ to $U$. A component of a vertex deletion is induced, and
the component property also persists. The drawing of $P$ lies in one face $F$ of $H'$. For
each simple cycle in $P$, the component $C_{U\cup A}$ lies on its large side. Since $H'$ is
connected and misses the cycle, its small disk misses $H'$. Thus the entire inside of the walk
belongs to $F$, and both designated faces merge into it. Every discarded component of $H-A$
was adjacent to $A$, so all discarded material except retained boundary vertices also lies in
$F$. All other faces of $H'$ are unchanged.

Trace the paths by keeping maximal nonempty segments in $H'$ and fixing the other edges. Each
new contact with the fixed edges is adjacent along a path to a vertex of $A$, so there are at
most $2|A|$ new contacts, all incident with $F$. A surviving vertex used by fixed edges alone
at a contact becomes a reserved singleton. Keep surviving reservations and record the order
and gains of segments. Reassign surviving sites of swallowed faces to $F$, and add one
boundary site if it has none. The number of faces falls and the number of incidences rises by
at most $2T(w)+1$.

If the sites $Z$ of a designated face $f$ are not free, a small side of order less than
$|Z|$ contains them. For $S\subseteq V(H)$ with $|S|<|Z|$, let $Q_S$ be the component
of $H-S$ containing $C_{U\cup S}$. Say that $S$ \defn{cuts off} $Z$ if $Q_S$ misses $Z$.
The separator of the small side cuts off $Z$. Choose a smallest such set $S$ and put
$r\coloneqq |S|$. Connectedness of $H$ gives $r\ge1$. Call a component of $H-S$ \defn{lost} if it contains a site of $Z$.

For $v\in S$, put $S^-\coloneqq S\setminus\{v\}$. By nesting of the large components,
$Q_{S^-}$ is the component containing $Q_S$ when $v$ is restored to $H-S$. Minimality
of $S$ implies that $Q_{S^-}$ meets $Z$. Thus $v$ has a neighbour in $Q_S$, and either
$v$ is a site or it has a neighbour in a lost component.

Delete the lost components to obtain $H'$. Every vertex of $S$ has a neighbour in $Q_S$,
and every component of $H-S$ has a neighbour in $S$. Thus $H'$ is connected and induced.
Each component of $J-(U\cup S)$ previously in $H$ lies in one component of $H-S$, and
$Q_S$ survives. Hence the component property of~(3) persists when we add $S$ to $U$. Each lost component,
together with the interiors of its edges to $S$, lies in the face $F$ of $H'$ containing
$f$, since it is connected, misses $H'$ and contains a boundary vertex of $f$. Every vertex of
$S$ lies on the boundary of $F$, either as a site of $f$ or through an edge to a lost
component. Replace
the old incidences of $Z$ at $f$ by those of $S$ at $F$. Reassign other surviving incidences
of swallowed faces to $F$. Trace the path segments as before, with contacts and singleton
reservations at $S$. The face count does not increase and the number of incidences decreases
by at least $|Z|-r$. After either operation all unchanged undesignated faces retain zero
quotient gain. Cellularity ensures that puncturing the designated faces suffices.

Using the bounds $z$ on initial faces and $w_0$ on initial incidences, let $w_{i+1}\coloneqq w_i+2T(w_i)+1$ for $i<z$
and put $w_*\coloneqq \max\{1,w_z\}$. The first operation occurs at most $z$ times. Between decreases in the number of faces the
second occurs at most $w_*$ times, since it lowers the number of incidences. Thus the total separator budget can be taken as
$u\coloneqq zT(w_*)+(z+1)(w_*+1)^2$. Set $\theta\coloneqq 1+T(w_*)>w_*$ and $c\coloneqq u+3\theta$. These
choices justify all orientations above. At termination $H-U$ contains $C_U$, so unbreakability
implies that at most $s$ vertices outside $U$ were discarded. Since $|U|\le u$, the induced
subgraph $H$ omits at most $s+u$ vertices in total. Every omitted edge has an end among these
vertices, since $H$ is induced. At every
traced contact a variable segment has a prescribed end or the vertex is reserved. The recorded
order therefore pastes any realizing residual family to the same fixed edges. At termination
distinct designated faces have disjoint site sets. Two face points incident with a common
vertex lie on a closed radial walk of length four, giving distance at most $2<T(w)$.
\end{proof}

A \defn{residual drawing modulo $L$} in a closed surface $\Sigma$ consists of the
retained subgraph $H$ of a drawn graph $J$ together with the data furnished for a
realization by Lemma~\ref{lem:genus-carrier-normalization}. In particular it records
the designated faces and their sites, the set $U$ of size at most $u$ and the tangle
$\mathcal T_H$ of order $\theta$. For $w$ final site incidences its \defn{required
distance} is $t\coloneqq T(w)$. Fix one point in each designated face.
Let $m,p$ be nonnegative integers with $p\le c$. A \defn{system of enclosures of
type $(m,p)$ for $H$} consists of two-port enclosures of facial cycles of $J$
and three-port enclosures of auxiliary cycles of $J$.
Their outer disks are pairwise disjoint and avoid $U$, the final sites and the
chosen face points. Each enclosure has at least $m$ surrounding cycles. For every
$X\subseteq V(J)$ of size at most $p$, each surrounding cycle disjoint from $X$
lies in $C_X$.

To \defn{suppress} an arc between consecutive ports is to replace it by a single
edge with the gain of the arc in that direction. A \defn{compression} of the system
is a subgraph of $H$ with the boundary arcs of its two-port and three-port cycles
suppressed. The designated faces and their sites are retained. The two-port cycles
bound facial digons and the three-port cycles bound facial triangles in the induced
drawing. The sites of these new faces are their ports.

For the next lemma put
\[
 \thmbound{lem:genus-well-compression}{t,u}\coloneqq (3t+1,u+3t).
\]
Here $u+3t\le c$ for the required distance $t$ of a residual drawing, since
$t\le T(w_*)<\theta$ and $c=u+3\theta$.

\begin{samepage}
\begin{lemma}
Let $\Sigma$ be a closed surface\label{lem:genus-well-compression} and let $L$ be a subgroup of a finite abelian group
$\Gamma$. Let $H$ be a residual drawing modulo $L$ in $\Sigma$, with separator bound
$u$, tangle order $\theta$ and required distance $t$. Every system of enclosures for
$H$ of type $\thmbound{lem:genus-well-compression}{t,u}$ has a cellular compression
whose labelling modulo $L$ extends flatly off punctures in the designated faces.
The designated faces and the new digons and triangles are distinct. The compression
admits a respectful tangle of order $t$ in which their site sets are free and their
pairwise distances are at least $t$.
\end{lemma}
\end{samepage}
\begin{proof}
Each outer perimeter lies in $C_U$, and its closed disk misses $U$. Since $J$ is
connected, every vertex in that disk is joined to the perimeter within the disk. Thus
the disk's entire graph lies in $C_U$ and survives in $H$. A designated face cannot
enter the disk, since its chosen point lies outside and its interior misses the perimeter.

In each disk keep the innermost perimeter and everything outside it. Inside it keep the
selected cycle and the portions of its spokes. Suppress the arcs of the cycle between
successive ports. The resulting graph $H'$ is connected, since any excursion across a
perimeter can be replaced by a perimeter arc and the retained spokes join the cycle to
the perimeter. The quotient transport restricts to $H'$ by giving each suppressed arc
its gain. The only punctures are still in the designated faces, away from the new
digons and triangles.

For $Y\subseteq V(H')$ with $|Y|\le3t$, choose in each disk its innermost perimeter $P$
missing $Y$. One exists among the $3t+1$ perimeters. All chosen perimeters lie in
$C_{U\cup Y}$, since $|U\cup Y|\le u+3t$. Excise from that component the open disks bounded by them, keeping their
boundaries, and call the remaining graph $B_Y$. It is connected by the same perimeter
detour. It is cellular because $C_{U\cup Y}$ is cellular, by the argument for $H$ in the proof of
Lemma~\ref{lem:genus-carrier-normalization}, and each excision replaces
the faces inside a perimeter by one disk face. Thus $B_Y$ is a nonempty cellular
subgraph of $H'-Y$. If $Y$ grows, the components shrink and the chosen perimeters
can only move outwards. Hence $B_{Y'}\subseteq B_Y$ whenever $Y\subseteq Y'$.

Orient each separation of order less than $t$ toward $B_Y$, where $Y$ is its separator.
For three such separations, the witness for the union of their separators misses all three
small sides. This gives a tangle. A short noose through $Y$ lies in a disk face of $B_Y$,
so its disk away from $B_Y$ is small. The tangle is respectful. Also $H'$ is cellular,
since it is connected and contains the cellular $B_\varnothing$.

A site set $Z$ is free if for every $|Y|<|Z|$ some vertex of $Z-Y$ can reach $B_Y$ in
$H'-Y$. For a designated face, freeness in $H$ gives a path from $Z-Y$ to $C_{U\cup Y}$
in $H-Y$. Otherwise the separation at $Y$ with every other component on its small side
would cover $Z$. Continue in $C_{U\cup Y}$ to $B_Y$ and stop at first contact. The path
starts outside all the selected disks and cannot enter a chosen perimeter before hitting
$B_Y$. Its retained prefix is therefore in $H'-Y$. For a new digon or triangle one of its
disjoint spokes misses $Y$ and joins a port to the chosen perimeter in $B_Y$.

Finally consider a closed radial walk of length less than $2t$ and let $Y$ be its primal
vertices. Its drawing misses each chosen perimeter because a radial edge can meet a
primal cycle only at a primal vertex. There are two cases.
If the walk is inside one chosen perimeter $P$, each simple cycle in its support bounds
a disk inside $P$. This disk is the small side since $B_Y$ is connected and contains $P$.
Thus the whole inside of the walk lies there and contains at most one of the specified
face points. If the walk is outside all chosen perimeters, it lies in the unchanged radial
drawing. Its cycles have the same insides as in $H$, since $B_Y\subseteq C_{U\cup Y}$
witnesses the same large side. Its inside therefore contains at most one designated face
point. For a new face point, its chosen perimeter and disk miss the walk.
Any small disk containing the point would therefore contain that perimeter in $B_Y$,
which is impossible. This proves all the distance bounds.
\end{proof}

\begin{lemma}\label{lem:rooted-frame-search}
Fix a closed surface $\Sigma$, a finite abelian group $\Gamma$, a subgroup $L\le\Gamma$,
and nonnegative integers $k,h,s$. For every nonnegative integer $w$ put
\[
 (j(w),T(w))\coloneqq \thmbound{cor:free-surface-frame}
   {\Sigma,\exp(\Gamma/L),w,w,w(|L|-1)} .
\]
Use $(w_*,u,\theta,c)\coloneqq \thmbound{lem:genus-carrier-normalization}{T,h+2k,h+2k}$.
Let $J$ be a $\Gamma$-labelled graph cellularly drawn in $\Sigma$ with at most $h$
defective faces modulo $L$, and let $\Pi$ be a pattern in $J$ of order at most $k$.
There is a polynomial-time algorithm that returns a realization of $\Pi$ or reports
failure. It succeeds if the hypotheses of Lemma~\ref{lem:genus-carrier-normalization}
hold at the fixed value $s$ and some realization has residual data from that lemma with
$w$ final site incidences and the following enclosures in $J$.
\begin{enumerate}
\item For each ordinary residual demand there is an ordered list of at most $|L|-1$
two-port enclosures of facial cycles whose gains belong to $L$ and whose subset sums
fill $L$.
\item There are $j(w)$ three-port auxiliary enclosures.
Together with the enclosures in the lists they form a system for $H$ of type
$\thmbound{lem:genus-well-compression}{T(w),u}$, where $H$ is the retained subgraph.
\end{enumerate}
When $L=\{0\}$ the lists in (1) are empty.
\end{lemma}
\begin{proof}
Enumerate a set of at most $s+u$ omitted vertices and let $H$ be the subgraph induced on
the remaining vertices. Keep connected cellular candidates. Enumerate visits to the omitted
vertices as in Lemma~\ref{near:apex}. This gives the ordered fixed edges, the residual
demands with their gains, and the singleton reservations where a path has no remaining edge.
If only reservations remain, test the guess directly.

Guess the designated faces of $H$ and at most $w_*$ site incidences, including the boundary
corner at each site. Each face needs a site, different faces must have disjoint sites, and
all names must be covered. In the parameter disk of each face draw a root disk touching
the specified sites in their chosen corners. For each ordinary residual demand choose at
most $|L|-1$ facial cycles of $J$ in order, with two distinct ports on each. Require them
to remain facial in $H$, their closed disks to be disjoint from one another and the root
disks, and each list's facial gains to lie in $L$ and have subset sums covering $L$.

For $w$ final incidences the size bound is
$\thmbound{prop:surface-relative-frame}{\Sigma,\exp(\Gamma/L),w,w,w(|L|-1)}$.
This is a constant for $w\le w_*$. Enumerate all graphs $F$ up to this bound with the
specified stars marked. Delete the star interiors and seek a subdivision of the exterior
graph with the sites and ports fixed. Delete the correction-cycle edges and their nonports
from the host, then guess the other branch images. Split every branch image into one
leaf for each incident host edge, keeping the edge attached to its leaf. For each
incidence of an exterior edge choose a distinct leaf at its branch image. Ordinary
disjoint paths between the chosen pairs, found by Robertson and Seymour~\cite{RSXIII},
give exactly the desired subdivisions after the leaves are identified again.
Any subdivision returned, together with the prescribed stars, embeds $F$. Enumerate
its disjoint routes for the given pairing that use the listed stars in order, avoid unused
root sites and open root disks, and avoid all other two-port disks. Expand the routes
through the subdivision, try either boundary arc at each visited star and paste the fixed
edges. Check the resulting paths and their gains in $J$. This test makes every output sound,
even for a guess that lacks the frame property.

For completeness choose the omitted vertices and visits from
Lemma~\ref{lem:genus-carrier-normalization}. Puncture each designated face in its root disk.
The residual paths are disjoint proper arcs with their true quotient gains after pushing
any incidental site contacts outwards in disjoint collars. With the lists in (1) and
the auxiliary enclosures in (2), Lemma~\ref{lem:genus-well-compression} gives a
$T(w)$-spread routing system of size at most $(w,w,w(|L|-1))$ in the compression.
By Lemma~\ref{cor:free-surface-frame}(1) it contains the exterior subdivision of a frame
with the property of Lemma~\ref{prop:surface-relative-frame}. Expanding the auxiliary
triangle edges lifts that subdivision to $H$ while avoiding the correction-cycle edges
and nonports. The search therefore succeeds for this abstract frame.

The returned subdivision also respects the marked disks and sites. The proof of
Lemma~\ref{cor:free-surface-frame}(2) applies to any such subdivision of the chosen
frame. It gives routes with the prescribed quotient gains and independent choices
on the correction cycles realizing the full gains. These routes pass the final check. Each
enumeration chooses boundedly many graph elements and each disjoint paths test has
fixed order.
\end{proof}

\section{The algorithm on a fixed surface}\label{sec:surface-algorithm}

\begin{theorem}\label{thm:all-abelian-surface-linkage}
Let $\Gamma$ be a fixed finite abelian group and let $g,k$ be fixed nonnegative integers.
Given a $\Gamma$-labelled graph $G$ of Euler genus at most $g$ and a reserved pattern
$\Pi$ in $G$ of order at most $k$, a deterministic algorithm finds a realization
of $\Pi$ or determines that none exists. For a graph of size $n$ its running time is
$O(n^{\thmbound{thm:all-abelian-surface-linkage}{\Gamma,g,k}})$ for an integer exponent
depending only on $\Gamma,g,k$.
\end{theorem}
\begin{proof}
Induct on Euler genus. On each of the finitely many closed surfaces $\Sigma$ of the current
Euler genus prove a stronger assertion for each $L\le\Gamma$, in increasing order of $|L|$.
The input has a drawing on $\Sigma$ whose labelling projects modulo $L$ to one extending
flatly off at most $h$ points. The demands are in $\Gamma$. The induction hypothesis holds for
every fixed $h$ and order, including new bounds needed for smaller subgroups. The case $L=0$
is included in the argument below. A valid pattern consisting only of reservations is
realized by those singleton vertices, so we may assume there is an ordinary demand.

At each induction step Lemma~\ref{lem:logic} extends the solver from unbreakable
structures to all inputs. Recursive linkage calls go only to earlier steps, whose
algorithms already work on all inputs, even for query parameters that depend on $s$.

Use Lemma~\ref{thm:linkage-rigid-skeletons} to reduce to rigid skeletons with patterns
of order at most $k_1$, where $(k_0,k_1)\coloneqq \thmbound{thm:linkage-rigid-skeletons}{k}$.
In the first case of Proposition~\ref{prop:surface-spqr-subgroup-guard}, enumerate
visits to the bounded deletion set by Lemma~\ref{near:apex} and use the lower-genus
algorithm on the remaining components. In the other case let $\varphi$ be the conjunction
of the bundle condition in that proposition and the linkage formula. Lemma~\ref{lem:logic}
reduces its evaluation to unbreakable structures. To test the bundle condition,
check the bundle marks and labels, the 3-connectivity of the reference graph and its
embeddability in $\Sigma$~\cite{MoharEmbedding}. In positive genus enumerate the
bounded deletion sets of Proposition~\ref{prop:surface-spqr-subgroup-guard} and
reject if any leaves only components of smaller Euler genus. In one embedding,
count the nonzero facial gains modulo $L$ and the adjacencies with incongruent
options. These tests take polynomial time. When they pass, the proposition gives
the cellular incidence drawing $J$. The unbreakability promise is on $J$ itself, and
its drawing has at most $H_0\coloneqq \thmbound{prop:surface-spqr-subgroup-guard}{\Gamma,k,h}$
defective faces modulo $L$.

This drawing has a bundle rendition. Enclose the subdivided options of each adjacency
in its edge strip, with the two ends as poles. The strips occupy disjoint sectors at each
vertex, and their torso is the simple 3-connected skeleton.

We choose the separation order before the size bound $s$ of Lemma~\ref{lem:logic}.
Assign each named vertex to one incident face and mark one site on each defective face
still without a name. Thus $z_0\coloneqq H_0+2k_1$ bounds both the number of designated faces
and the number of site incidences. Use the functions $j,T$ of
Lemma~\ref{lem:rooted-frame-search} for this surface, group and subgroup, and set
\[
 (w_*,u,\theta,c_0)\coloneqq \thmbound{lem:genus-carrier-normalization}{T,z_0,z_0},
 \qquad (m,p)\coloneqq \thmbound{lem:genus-well-compression}{T(w_*),u}.
\]
At most $u+w_*+z_0$ disjoint outer disks can meet the separator set, the final sites
or one chosen point in each designated face. After these losses we need at most
$w_*(|L|-1)$ copies of any gain class and $j(w_*)$ auxiliary enclosures, since $j$
is nondecreasing. Hence in Lemma~\ref{lem:mixed-represented-digons} take
\[
 q\coloneqq u+w_*+z_0+\max\{w_*(|L|-1),j(w_*)\}.
\]
In particular $q>H_0$. Write $(b,d,N)\coloneqq \thmbound{lem:mixed-represented-digons}{\Sigma,\Gamma,m,p,q,s}$.
The bound $b$, and on the sphere also $d$, is independent of $s$.
Choose $c\ge c_0$, with $c\ge2d+1$ on the sphere, and then set
$s\coloneqq \thmbound{lem:logic}{\varphi,c}$.

It now suffices to solve the formula when $J$ is $(s,c)$-unbreakable, for a fixed $s$. Graphs
with at most $3s+c$ vertices are tested directly. On the sphere use the bounded-treewidth
algorithm up to the threshold $N$. At positive genus enumerate bounded deletions up
to $\max\{s+c,N\}$. If a deletion leaves all components of smaller Euler
genus, use the outer induction as before. Otherwise the drawing has the required
representativity by cutting and capping. Its facial boundaries are simple by
Proposition~\ref{prop:surface-spqr-subgroup-guard}.

Let $K$ be the subgroup generated by the classes in the set $A$ of
Lemma~\ref{lem:mixed-represented-digons}. Each class occurs on $q>H_0$ distinct
faces, so $K\le L$. We try both possibilities without computing $A$.

First run the search of Lemma~\ref{lem:rooted-frame-search} with parameters $k_1,H_0$.
It is complete when $K=L$. Indeed, normalize any realization by Lemma~\ref{lem:genus-carrier-normalization}.
For each ordinary residual demand choose a list of at most $|L|-1$ gains
from the classes in $A$, allowing repetitions, whose subset sums fill $L$, as before Lemma~\ref{cor:free-surface-frame}. The choice of $q$ supplies
these enclosures and the auxiliary enclosures outside all forbidden points.
Since $T(w)\le T(w_*)$, our choices of $m,p$ give the required type. For arbitrary
$K$ the search still verifies every linkage it returns.

In the second search put $(r,v)\coloneqq \thmbound{lem:confined-face-regions}{s,d}$.
Enumerate all proper subgroups $K'<L$ and at most $b$ original faces. Let $X$ consist
of the primal vertices within radial distance $r$ of the face centres. Test every
component of $J-X$. On the sphere accept a component of treewidth at most $v$. At
positive genus accept a planar component. In either case also accept a component whose given
subdrawing is flat modulo $K'$ off at most $b$ points, by
Lemma~\ref{thm:find-flatness-punctures}. If all components pass, use
Lemma~\ref{cor:mixed-radial-ball-peeling}, whose subdrawing queries have order at most
$\thmbound{cor:mixed-radial-ball-peeling}{\Sigma,\Gamma,k_1,b,r}$. On subdrawings avoiding $X$ the accepted conditions
persist. In the width case use bounded-treewidth linkage. In the planar case use the already
completed Euler-genus-zero algorithm for the full group. In the flatness case use the inner
induction at $K'$ on the same surface. Thus this branch gives an exact decision. When $K<L$,
Lemma~\ref{lem:confined-face-regions} guarantees that the choice of $K'=K$ and of the faces
of that lemma is accepted. For
$L=0$ we always have $K=L$, so the first search alone is complete and no smaller-subgroup call
is needed. If neither search succeeds, the input is negative.

All steps are deterministic and polynomial for the fixed parameters. Lemma~\ref{lem:logic} now
removes the unbreakability assumption and the block tables give the decision in the input
graph. Taking $L=\Gamma$ removes any condition on its labelling. An embedding in a specified
fixed surface can be found in polynomial time by the algorithm of
Mohar~\cite{MoharEmbedding}. Try the finitely many
surfaces of Euler genus at most $g$. To recover a witness after a positive decision,
test each edge once and delete it if a linkage remains. Every surviving edge belongs
to every remaining linkage, so the components containing edges are precisely the
ordinary paths. Restore the reserved vertices as singletons.
\end{proof}

\section{Passage to arbitrary graphs}\label{sec:journal-global}

For a subgroup $\Lambda\le\Gamma$ and a positive integer $n$, the \defn{labelled clique}
$K(\Lambda,n)$ has vertex set $[n]$ and an edge $(i,j,a)$ with tail $i$, head $j$ and
label $a$ for each pair of distinct indices $i,j\in[n]$ and each $a\in\Lambda$.
A \defn{complete $\Lambda$-model of order $n$}
in a $\Gamma$-labelled graph is a model of this clique by disjoint \defn{branch trees}
$T_1,\ldots,T_n$ and distinct \defn{representative edges}. The representative of $(i,j,a)$
has its tail in $T_i$ and its head in $T_j$. One shift makes the tree edges zero and
gives every representative its specified label. A connected component outside a block has at
most one neighbour in that block, since a path between two different neighbours would extend
the block by an ear. In particular, a simple path between its vertices stays in the block. For
$n\ge3$, all representative edges of a complete model lie in a single block. Any three
distinct tree indices give a cycle through any three chosen representatives on the three
pairs, by adjoining tree paths. These triangles connect all representative edges into one
block. Pruning each tree to the minimal subtree through its representative ends puts the
entire model in that block.

\subsection{A complete minor and a small separator}
For a nonnegative integer $m$, put $\thmbound{lem:journal-clique-interface}{m}\coloneqq 3m+3$.

\begin{lemma}\label{lem:journal-clique-interface}
Fix a finite abelian group $\Gamma$, a subgroup $\Lambda\le\Gamma$, and nonnegative integers
$m,n$ with $n\ge\thmbound{lem:journal-clique-interface}{m}$. Let $C$ be a $\Lambda$-labelled graph supplied with a complete
$\Lambda$-model of order $n$, and let $W\subseteq V(C)$ have size $m$. One can find in
polynomial time a separation $(A,B)$ of $C$ with separator $S\coloneqq V(A)\cap V(B)$ of order
at most $m$ and $W\subseteq V(A)$. All model trees disjoint from $S$ lie in $B-A$.
Every pattern of ordinary demands on distinct vertices of $S$ with prescribed gains in $\Lambda$
has a realization in $B$ avoiding all other vertices of $S$.
\end{lemma}
\begin{proof}
For $|Q|<n$, all model trees missed by $Q$ belong to one component $U_Q$ of
$C-Q$, since there is a representative edge between every two of them. Call a side
with separator $Q$ \defn{small} if it misses $U_Q$. Choose a minimum-order small
side $A$ containing $W$, with opposite side $B$ and separator $S$ of size $p$.
Taking $Q=W$ shows that $p\le m$. For each candidate separator $Q$ of
order at most $m$, it suffices to test whether $V(C)\setminus U_Q$ contains $W$.
This is the largest small side with that separator.
In $B$ the trees missed by $S$ give a model of order at least $n-p>2p$.
No small side of order less than $p$ in $B$ contains $S$, since adjoining $A$
to it would give a smaller-order small side in $C$ containing $W$.

It suffices, then, to prove the following routing claim. Suppose a
$\Lambda$-labelled graph contains a complete $\Lambda$-model of order $N>2p$,
and a set $Y$ of $p$ vertices lies in no small side of order less than $p$.
Every pairing of vertices of $Y$ with gains in $\Lambda$ can be routed while
avoiding the unused vertices of $Y$. We may assume $p>0$.

Among the small vertex sides of order $p$ containing $Y$ there is a largest, say
$D$. Indeed, $Y$ itself is one such side. The union and intersection corners of
two such sides both contain $Y$ and miss a common whole model tree, since their
separators together meet at most $2p<N$ trees. Both corners are therefore small.
Their orders sum to $2p$ and neither is below $p$, so the family is closed under unions.

Take a counterexample with the fewest model-tree edges. If a tree edge $uv$ has
$v\notin D$, shift at $v$ to make the edge zero and contract it. Since $Y\subseteq D$,
the shift is zero on the terminals and they remain distinct. All labels stay in $\Lambda$. A small side of
order less than $p$ containing their images would lift to one containing $Y$ by
restoring both ends of $uv$ wherever its contracted vertex occurs. A whole model
tree missed by the lifted separator witnesses smallness. The lifted order is at
most $p$. Equality requires the contracted vertex in the separator, putting $v$
on a small side of order $p$ containing $Y$, contrary to $v\notin D$. The
contracted instance therefore has a routing by minimality. Restoring the zero
edge on the path through the contracted vertex, if needed, lifts that routing.

Consequently the model trees missed by the separator of $D$ are singleton
vertices outside $D$. Choose $p$ of them as $Z$. For every $|Q|<p$, the
component $U_Q$ meets $Z$ because some chosen tree survives, and meets $Y$
because the small side outside $U_Q$ cannot contain $Y$. Menger's theorem
gives $p$ disjoint $Y$--$Z$ paths. Between two vertices of $Z$ the representative
edges in each direction realize all of $\Lambda$. Their labels belong to
$\Lambda$ and a common translation gives all the prescribed model labels.
For each demand join its two paths by the edge giving the desired total gain.
Discard the paths from unused vertices of $Y$. This contradicts the choice of a
counterexample and proves the claim. Applying the claim to $Y=S$ in $B$ proves the lemma.
\end{proof}

\begin{lemma}\label{thm:ordinary-exclusion-suffices}
Let $\Gamma$ be a finite abelian group and let $k$ be a nonnegative integer. There is a
positive integer $h\coloneqq \thmbound{thm:ordinary-exclusion-suffices}{\Gamma,k}$ such that
a polynomial-time linkage decision algorithm for patterns
of order at most $k$ in $\Gamma$-labelled graphs excluding the ordinary complete minor
$K_h$ suffices for a polynomial-time decision algorithm for all patterns of that order
bound in arbitrary $\Gamma$-labelled graphs.
\end{lemma}
\begin{proof}
Write $u\coloneqq 2k$, $l\coloneqq \lfloor\log_2|\Gamma|\rfloor$, and choose $N_l\coloneqq 8(u+1)$. For
$j=l-1,\ldots,0$ set $t_j\coloneqq 8|\Gamma|^2N_{j+1}$ and $N_j\coloneqq 8(t_j+u+1)$. Let
$a_*\coloneqq u+\max(\{0\}\cup\{t_j:j<l\})$ and $c\coloneqq 3a_*+1$. By Lemma~\ref{lem:logic} it suffices to
solve decoded instances in which every separation of order at most $c$ has a strict side of
size at most $s\coloneqq \thmbound{lem:logic}{\varphi,c}$, where $\varphi$ is the fixed
formula for the demand type and target tuple under consideration. Choose $M\ge\max(N_0,s+c+3)$. By Geelen and
Gerards~\cite[Theorems~1.3 and 1.4]{GeelenGerards}, an ordinary clique minor of sufficiently
large order $h$ supplies a complete $\{0\}$-model of order $M$. Their construction in
Theorem~1.4 normalizes the branch trees, finds a monochromatic ordered subclique by
enumeration, shifts one side of its complete bipartite subgraph and contracts a matching.
These operations take polynomial time for fixed parameters. The starting ordinary clique model
or a certificate of exclusion is found by the fixed-minor testing algorithm of Robertson
and Seymour~\cite{RSXIII}. In the latter
case use the assumed solver. In the former case retain the representative edge set $E_0$ of
the entire neutral model.

Begin with $N_0$ of these trees and subgroup $\Lambda\coloneqq \{0\}$. At stage $j$ we have a complete
$\Lambda$-model of order $N_j$ with representatives in $E_0$. If $\Lambda\ne\Gamma$, apply the
cited Theorem~1.3 to $4t_j$ of its trees. The initial minor and each subsequent minor
have their representative edges oriented as in the paper of Geelen and
Gerards~\cite{GeelenGerards}. In polynomial time this gives either a complete
model of order $N_{j+1}$ over a strictly larger subgroup with representatives among the
previous ones, or a set $Z$ of fewer than $t_j$ vertices for which the block of $G-Z$
containing more than half the sampled representatives shifts into $\Lambda$. In the former
case continue. Each step at least doubles the subgroup size. In the latter case every
surviving tree of the full model can be pruned into this block. Indeed all edges between
untouched trees lie in one block, and more than $3t_j$ of the sampled trees are untouched. The
representatives between them constitute more than half the sample because
$\binom{3t_j}{2}>\tfrac12\binom{4t_j}{2}$ for $t_j\ge 2$. On reaching $\Gamma$ stop with
$Z\coloneqq\varnothing$.

Delete $X$, the union of $Z$ with the named vertices. Put $a\coloneqq |X|\le a_*$. At least $n\coloneqq N_j-a$
trees remain in a single block $C$ of $G-X$, and $N_j\ge8(a+1)$. This block still shifts into
$\Lambda$ because a block containing the surviving representatives after further deletions
lies in the old block. Shift $C$ to $\Lambda$ by normalizing a spanning tree, extending the
potential by zero outside. None of the named vertices changes its potential.

For any set $D$ of $d\le c$ vertices, the original trees missed by $D$ are
pairwise adjacent in $G-D$. Their component $Q_D$ has at least $M-d>s$
vertices, so unbreakability bounds the rest of $G-D$ together by $s$ vertices.
Each other component meets only the original trees hit by $D$, so it contains
at most $d(d-1)$ edges of $E_0$, one per ordered pair of distinct trees.
Consequently any surviving complete model of order $r>d$, with $r\ge2$ and
representatives in $E_0$, lies in $Q_D$, since it is connected and uses at least
$r(r-1)$ such edges.

Now guess all edges of a solution incident with $X$. There are at most $2a$.
Keep $C$, $X$ and the components of $(G-X)-V(C)$ reached by a guessed edge.
Delete all unguessed edges incident with $X$. Each exterior component has at most one attachment to $C$.
Since the terminals lie in $X$, a solution visiting the component must use an edge
from it to $X$ or pass twice through its attachment. The latter is impossible on a path.
Thus some guess preserves a solution. Let $D$ consist of $X$ and the attachment of
such a component, if any. At least $N_j-a-1>|D|+1$ of our model trees in $C$ miss $D$, so
$Q_D$ meets $C-D$. The exterior component is a different component of $G-D$
and thus has at most $s$ vertices. The retained exterior has
at most $a+2as$ vertices outside $C$ and meets $C$ in a set $W$ of at most
$2a$ vertices, one per guessed edge to $C$ or per retained component.

Since $N_j\ge8(a+1)$ and $|W|\le2a$, the surviving model has order at least
$N_j-a\ge\thmbound{lem:journal-clique-interface}{|W|}$. Apply
Lemma~\ref{lem:journal-clique-interface} in $C$ with this $W$, and let $S\coloneqq V(A)\cap V(B)$. Put $D\coloneqq X\cup S$ and $d\coloneqq |D|\le3a$. At least $N_j-a-|S|>d+1$ of our model trees
survive strictly on the $B$ side, with representatives in $E_0$, so they lie in $Q_D$.
Paths of $G-X$ between vertices of $C$ stay in $C$, and no path of $C-S$ joins $A-S$ to
$B-S$. Hence $A-S$ lies outside $Q_D$, giving $|V(A)\setminus S|\le s$. Replace $B$ by one
edge for each pair of $S$ and each possible traversal gain in $\Lambda$. Each excursion in $B$
shortens to such an edge. Conversely each maximal run of replacement edges has its total gain
in $\Lambda$. Realize these runs in $B$ by the preceding lemma avoiding the other roots, and
join them to the exterior paths. This preserves disjointness and the required gains. The
replacement graph has boundedly many vertices, so a direct enumeration of simple vertex
sequences and available labels decides the pattern. All choices made have polynomial cost for
the fixed parameters. Applying Lemma~\ref{lem:logic} completes the reduction. The finitely
many demand types and tuples for the fixed group and order can use the largest of their clique
thresholds.
\end{proof}

\subsection{The reduction to a surface}
\begin{theorem}\label{thm:bounded-genus-host-suffices}
Let $\Gamma$ be a fixed finite abelian
group and let $k$ be a fixed nonnegative integer. There are integers
$(g_*,k_*)\coloneqq \thmbound{thm:bounded-genus-host-suffices}{\Gamma,k}$ with the
following property. A polynomial-time linkage decision algorithm for patterns
of order at most $k_*$ in $\Gamma$-labelled graphs supplied with a drawing
in a closed surface of Euler genus at most $g_*$ suffices for a polynomial-time
decision algorithm for patterns of order at most $k$ in arbitrary
$\Gamma$-labelled graphs.
\end{theorem}
\begin{proof}
By Lemma~\ref{thm:ordinary-exclusion-suffices} it suffices to treat graphs
excluding $K_h$, where $h\coloneqq\max\{5,\thmbound{thm:ordinary-exclusion-suffices}{\Gamma,k}\}$.
Delete all loops and subdivide each remaining edge on its own new vertex,
orienting the two halves from the old tail towards the old head.
Give the first half the original label and the second half label zero.
The resulting graph is simple and its paths between old vertices
correspond by suppression to the original paths with the same gains.
It also excludes $K_h$. Indeed every branch set of a $K_h$ model in
the subdivision must contain an old vertex, since a new vertex has
degree two and distinct new vertices are nonadjacent. Suppressing
the subdivision vertices then gives a model on disjoint connected
sets of old vertices in the original graph.

Choose the cut order $K$ specified for $h$ in
Appendix~\ref{app:construct}. For each fixed type of pattern of order
at most $k$ and each tuple of gains, let $\psi$ be its formula from
Lemma~\ref{lem:logic} conjoined with simplicity and exclusion of
$K_h$ in the decoded graph. The latter is expressible by excluding
$h$ disjoint nonempty connected vertex sets adjacent in every pair.
By that lemma it suffices to decide $\psi$ on the promised structures.
Reject malformed encodings, nonsimple decoded graphs and decoded
graphs with a $K_h$ minor. On each remaining input every separation
of order at most $K$ in the decoded graph has a strict side of size
at most $s\coloneqq \thmbound{lem:logic}{\psi,K}$.
Appendix~\ref{app:construct} gives a polynomial-time reduction on
these graphs to the stated surface queries, with the genus and
orders bounded in terms of $h,s,\Gamma,k$. Take coordinatewise
maxima over the finitely many demand types and gain tuples to obtain the claimed bounds.
\end{proof}

\begin{proof}[Proof of Theorem~\ref{thm:main}]
Theorems~\ref{thm:bounded-genus-host-suffices} and~\ref{thm:all-abelian-surface-linkage}
together decide whether $G$ has a $\Pi$-linkage. After a positive answer, the edge
deletions in the proof of Theorem~\ref{thm:all-abelian-surface-linkage} recover the paths
with at most $|E(G)|$ further calls. The resulting polynomial has an exponent depending
only on $\Gamma$ and $k$.
\end{proof}

\section*{AI Declaration}
This proof was found by ChatGPT 6.2 Sol, ChatGPT 6.2 Astra, and Claude Fable 5.2.  The process began with asking the HAL 9000\footnote{HAL 9000:=$\{$ChatGPT 6.2 Sol, ChatGPT 6.2 Astra, Claude Fable 5.2$\}$.} to read the author's PhD thesis, and took about two weeks in total, with guidance from the author.  However, the key ideas of using unbreakable tree-decompositions~\cite{LRSZ} and the coarse Menger theorem for bounded genus graphs~\cite{BPPcoarseMenger} are entirely due to the HAL 9000.

\section*{Acknowledgements}
The author is extremely grateful to Jim Geelen for his supervision and for several key ideas that led to this work. He also thanks Guoli Ding, Bert Gerards, Thor Johnson, Stephan Kreutzer, Neil Robertson and Paul
Seymour for their contributions to the original work, and Bertrand Guenin, Naomi Nishimura,
Bruce Richter and Xingxing Yu for their careful reading and suggestions.

\appendix
\section{Constructing the surface reduction}\label{app:construct}
Fix a finite abelian group $\Gamma$, integers $h\ge5$ and $k\ge0$, and a simple
$\Gamma$-labelled graph $G$ excluding $K_h$. We construct the reduction used
in Theorem~\ref{thm:bounded-genus-host-suffices} for patterns of order at most $k$.
All wall and embedding computations below use the underlying graph.

\subsection{Obtaining a drawing and disjoint enclosing disks}
The following near-embedding terminology is that of
Fiorini, Joret, Weltge and Yuditsky~\cite{FJWY}. After deletion of
a set of vertices called \defn{apices}, the graph consists of a drawn graph $F$ and edge-disjoint
pieces. Each piece has a disk whose interior misses $F$. It meets the
rest of the graph only in its \defn{society}, the vertices of $F$
on the boundary of its disk, taken in boundary order. The disks have
disjoint interiors. The pieces are of two kinds. A \defn{small piece} has at most three society
vertices, called its \defn{ports}. A \defn{vortex} has $m\ge1$ society vertices
and a path decomposition $(U_1,\ldots,U_m)$ with
its $i$th society vertex in $U_i$. A vortex may itself have at most three
society vertices. Its \defn{adhesion} is the
maximum size of $U_i\cap U_{i+1}$.
Different pieces may share society vertices. The \defn{augmented drawing}
$F'$ joins the ports of each small piece into a clique, drawing any
required new edges in the corresponding disks. A small piece is
\defn{properly attached} if any two distinct ports can be joined inside
it by a path avoiding the other ports, and for any three distinct
ports $x,y,z$ such an $x$--$y$ path and such a $y$--$z$ path can be
chosen to meet only at $y$.

We use Theorems~15 and~19 in the preprint numbering
of~\cite{FJWY}. Their first theorem, adapted there from Kawarabayashi, Thomas and
Wollan~\cite[Theorem~2.11]{KTW} in the terminology of~\cite{DKMW}, gives initial bounds $a_0,u_0,d_0,g$ on the numbers of
apices and vortices, their adhesions and the Euler genus, respectively.
These depend only on $h$. Their second theorem gives the bounds
\[
       u\coloneqq u_0+g,\qquad a\coloneqq d_0+u_0+g
\]
on the final number of vortices and their adhesions. For each
$\beta\ge3$ it gives an apex bound $A_\beta$ depending only on $h$
and $\beta$. In particular these three bounds and $g$ are independent
of the requested height $r\ge\beta$ of the final wall. For fixed parameters both constructions
take polynomial time. Fiorini et al. state the refinement as implicit in the proof of
Diestel et al. and specify the stopping point giving these bounds~\cite[Theorem~19 and
footnote~3]{FJWY}.

The \defn{elementary wall} of height $t\ge2$ is obtained from the
grid on $\{1,\ldots,2t\}\times\{1,\ldots,t\}$ by deleting the
vertical edge $(i,j)(i,j+1)$ whenever $i+j$ is even and then
discarding the two vertices of degree one. A \defn{wall} of height
$t$ is a subdivision of this graph. Its \defn{perimeter} is the
subdivision of the exterior facial cycle in the grid drawing.
An embedded wall is
\defn{flat} if its perimeter bounds a disk containing the wall.
Take $b\coloneqq 2a+2$, $\beta\coloneqq b+6$,
and choose the cut order
\[
 K\coloneqq \max\{a_0+3,\ A_\beta+u(b+1)+2b\}.
\]
This depends only on $h$. Once an unbreakability size bound $s$ is given, take
$r\ge\beta$ large enough that a flat wall of height $r$ has more
than $s+3$ vertices strictly inside its perimeter. We now work on $(s,K)$-unbreakable inputs and request this height in the
final drawing.

We first prepare a flat wall for Theorem~19 of~\cite{FJWY}. Let $t$ be
large enough that any graph of treewidth above $t$ drawn in a surface
of Euler genus at most $g$ contains a flat wall of the input height
required there for $\beta,r$ and our initial bounds. Indeed the planar
minor theorem~\cite{RSV} gives any fixed wall as a minor. Since the wall
is subcubic, its minor model contains a subdivision. Lemma~7 of~\cite{DKMW} gives the
required flat subwall in a sufficiently large wall. Algorithmically, find the
minor by~\cite{RSXIII}, prune its branch trees, and enumerate subsets of the edges of the fixed elementary wall,
retaining their entire subdividing paths. Test the resulting subwalls and the
disk sides of their perimeters in the drawing.

Apply Theorem~15 of~\cite{FJWY} with excluded minor $K_h$, requesting a
wall $W'$ of treewidth greater than $\max\{t,s+3\}$. A sufficiently large
height suffices since walls contain arbitrarily large grid minors. The
input wall of the requisite height is again found by~\cite{RSV,RSXIII},
with smaller treewidth handled by Lemma~\ref{prop:width}.
Since $G$ excludes $K_h$, the theorem gives $W'$ in the augmentation. Its virtual edges
can be replaced by paths in the small pieces to obtain a subdivision of
$W'$ in the union of $F$ and the small pieces~\cite[Theorem~15 and Section~4.3]{FJWY}.
In particular $F$ has more than $s+3$ vertices.

Refine the small pieces to be properly attached as in
Lemma~3 of~\cite{DKMW}. To see that this takes polynomial time, consider
a piece that fails proper attachment. If every component minus the ports
has at most two neighbouring ports, draw the edges between ports and
give each component its own disk with its neighbouring ports. Two-port
disks run in parallel strips along the appropriate side of the old disk.
For a single neighbouring port use a sector at that port. A component missing all
ports gets a disk in the interior.
The resulting pieces are properly attached. Otherwise one component has
neighbours at all three ports, so each pair has a path avoiding the third.
Failure now means that some port $y$ has no two disjoint paths, apart from
$y$, to the other two. By Menger's theorem a vertex $v\ne y$ separates
$y$ from them. The paths avoiding the third port force $v$ to be private.
Draw $v$ in the disk and split off the component at $y$ of
the piece minus $v$ as a two-port piece on $y,v$. The remainder has the
other two ports and $v$ as its ports. Each such split promotes a vertex,
so component and cutvertex searches suffice. The bounds $a_0,u_0,d_0,g$
are preserved.

Let $H$ be the union of the drawn graph and its small pieces. The
refinements merely repartition $H$, so it still contains the lifted wall.
The ports of each new small piece, together with the apices, separate
its private vertices from more than $s$ drawn vertices. Hence the piece has size at most $s+3$.
In the new augmentation $F'$ its ports form a clique. Thus a tree
decomposition of $F'$ extends to one of $H$ by attaching each small piece
as a leaf bag at a bag containing its ports and ignoring virtual edges.
Writing $\operatorname{tw}$ for treewidth, this gives
\[
  \operatorname{tw}(H)\le\max\{\operatorname{tw}(F'),s+2\}.
\]
The lifted wall forces $\operatorname{tw}(H)>\max\{t,s+3\}$.
We can therefore find the required flat wall directly in $F'$ and
apply Theorem~19 of~\cite{FJWY}.

Write $F'$ now for the augmented drawing of the output. By the definition of a good
near-embedding in~\cite[Section~5.1]{FJWY} it contains a flat
wall of height $r$. Around each vortex there exist $\beta$ disjoint cycles in $F'$ bounding nested disks
that contain its society. The outer disks are pairwise disjoint and
miss the disk of the wall. All small pieces are properly attached.
Every small piece, and every vortex with at most three society vertices, has at most $s$
private vertices. Indeed its society and the apices together have size at most $A_\beta+3\le K$.
They separate its private vertices from more than $s$ wall vertices, since the wall
has more than $s+3$ vertices strictly inside its perimeter, all in $F$, and at most three
belong to the society.
We recover the opening arcs from the augmented drawing. The promised enclosing cycles
and exterior wall certify that this search succeeds.

\subsection{Finding the opening arcs}
Set aside the at most $u$ vortices with at most three society vertices.
We keep these vortices in full as compulsory exceptions below. For every remaining vortex draw an
auxiliary polygon through its society in order, just inside the old disk.
Use the smaller polygonal disk as its new disk and mark its open face by $z$.

Try one exterior point $z_*$ per face of the map with all polygons inserted.
For each vortex, breadth-first search from its own $z$ assigns each face $f$
a distance $d(f)$, allowing a step between faces sharing a vertex.
At level $1\le i\le b$, take the edges separating faces of distance at
most $i$ from the other faces. They form an even subgraph, since the
changes between the two kinds of face occur in pairs around a vertex.
Decompose it into simple cycles. Choose any cycle of $F'$ from this
decomposition bounding a disk containing the closed polygonal disk and
missing $z_*$. Discard the trial if a choice fails. Check that the disks
are strictly nested in level order for each vortex, the outermost disks
are pairwise disjoint and miss the societies set aside, and more than
$s+3$ vertices of $F'$ lie outside their union. These tests are polynomial.
For the disk test first check that the cycle is two-sided and separating.
After cutting along it, a component is a disk exactly when it has one
boundary circle and Euler characteristic one.

A successful trial gives the opening arcs as follows. The incident face distances at a vertex differ by at most one. At either end of
a level-$i$ edge their minimum is therefore $i$.
Along a shortest face route from $z$ to $z_*$ the contact vertices thus
have successive minimum incident distances $0,1,\ldots$. Such a route
can be drawn as a simple \defn{normal arc}, one whose only graph contacts
are vertices, since its faces and contact vertices are distinct. It must
cross each selected cycle, and can do so only at the contact of that level.
The subarc $L$ from its first contact, a society vertex, through the level-$b$
contact therefore meets each cycle once and has no other graph contacts.

To see that a trial succeeds, take $z_*$ inside the promised wall perimeter.
Of the $\beta$ enclosing cycles of a vortex, the outer $\beta-1$
strictly enclose its society and hence its polygon. Reaching the outermost
cycle from $z$ along a normal arc requires at least $\beta-1$ contacts.
The closures of faces of distance at most $b$ therefore lie strictly inside
the outermost disk. Indeed one may approach any point of such a closure
within the face after at most $b$ contacts and meet at most one further
vertex at the point itself. Thus the level edges lie there and miss
all other polygons. They also avoid their own polygon, whose vertices have
minimum incident distance zero. A generic curve from $z$ to $z_*$ crosses
the level edges an odd number of times, since its end faces have opposite
markings. Some cycle of the decomposition therefore bounds a disk containing
the polygon and missing $z_*$. Any such disk is inside the promised outer
disk. Different levels have disjoint cycles, and the shortest arc crosses
them in level order, so their disks are strictly nested. Distinct vortices
have disjoint promised outer disks, and the wall supplies the required
exterior vertices. The successful trial need not identify this wall.

\subsection{Preserving the gains at the small pieces}
Replace each small piece by a clique on its ports. With two ports give its edge
the gain of a joining path. With three ports there are two possibilities.
If a component off the ports sees all three, take a tree with the ports as its leaves
and give each edge of the clique the gain of the corresponding tree path.
Otherwise proper attachment gives a path between each pair avoiding the third port.
These paths form a subdivided triangle, since a component used by two pairs would see all three ports. Use its edge gains. In all cases gains are taken in the
chosen edge directions. Retain all drawn edges and distinguish edges from different
pieces. Draw the resulting multigraph $F^+$ along $F'$, with one copy on each old
trace and the others in narrow parallel strips. The opening arcs keep their vertex
contacts and the selected cycles keep their drawings.

Apply Lemma~\ref{near:apex} to the final apex set $A$, combining the residual demands
in $G-A$ into one pattern of order at most $\thmbound{near:apex}{k,|A|}$.
Mark every piece with a named private vertex and keep it as an exception. In every other piece
a linkage has at most one maximal run of edges, since two runs would require four
distinct ports. Shorten each run to the edge between its end ports. To recover the
original gain of a path, list in order the gain lost at each shortening. Whenever
two partial sums, including the initial zero, agree, erase the intervening terms. This preserves the sum and
leaves at most $|\Gamma|-1$ terms. Keep their pieces as further exceptions and use
the original runs in them. Each path now has its prescribed gain, and shortening
elsewhere preserves disjointness.

Conversely, a run through a replacement clique lifts to the chosen path or subdivided
triangle. In the tree case the gains along a run telescope to the gain of the tree path
between its ends. This path avoids the third port, which is a leaf. The lifted runs
use only private vertices in addition to their original ports, so they preserve
disjointness. Thus it suffices to retain the marked pieces and at most $|\Gamma|-1$
further pieces per ordinary residual demand.

Enumerate all sets of exceptions of the indicated size that contain every marked piece,
adding every vortex
previously set aside to every set. In each
branch restore those pieces, discarding their replacement edges,
and enumerate visits to all their vertices, including their ports,
by Lemma~\ref{near:apex}. Their bounded sizes make this a
bounded-order operation. In the resulting residual tests use the
intact drawing $F^+$ and reserve all exceptional ports by singleton demands. The private
vertices of the exceptions are gone, and every restored replacement
edge for them has both ends forbidden. In particular their edges
remain available in the drawing for structural arguments.

\subsection{Rooting the attachments}
Add the drawing vertices on each opening arc $L$ to the new apex set. Order the
surviving society from the cut at the society end of $L$. At each surviving society
vertex $w_t$, at most $2a+1$ vertices, including $w_t$, separate the strictly earlier
from the strictly later society vertices in the old vortex. For the first or last
vertex there is nothing to separate. Otherwise the earlier sites form a cyclic interval.
In the linear order of a promised vortex decomposition, either this interval or its
complement is consecutive. Delete
the intersections of consecutive bags at its two ends, taking the empty set at an
end of the decomposition. Each intersection has size at most $a$. Every surviving
vertex occurs on only one side, since its bags are consecutive, and no edge joins the
two sides, since its ends share a bag. Adding $w_t$ gives the claimed separator.
Only its existence is needed.

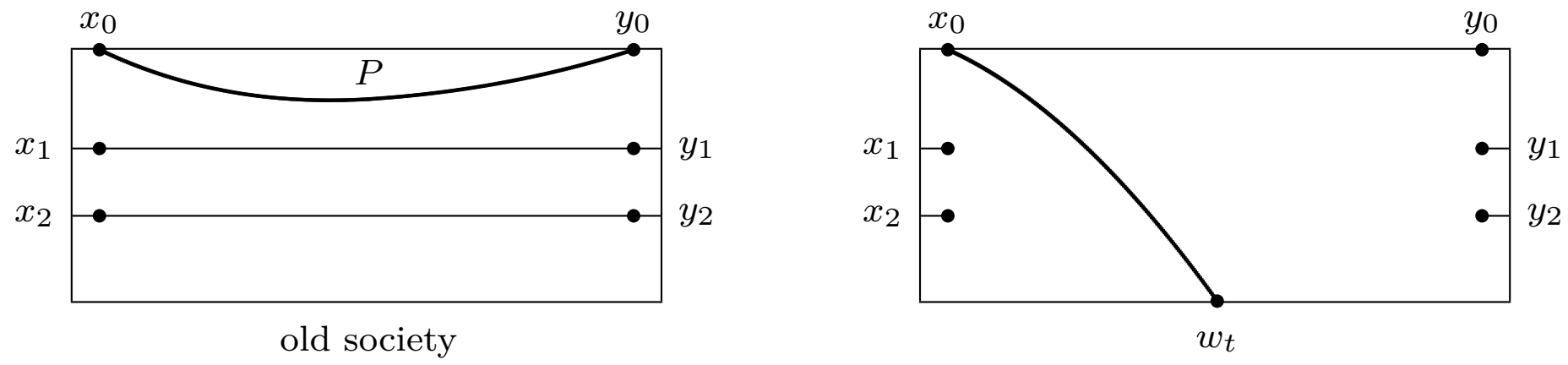
\begin{figure}[tbp]
\centering
\setlength{\unitlength}{1pt}
\scriptsize
\newcommand{\cutpanel}{%
  \thinlines
  \put(19,15){\framebox(131,56){}}
  \put(25,71){\circle*{2.8}}\put(144,71){\circle*{2.8}}
  \put(25,77){\makebox(0,0){$x_0$}}\put(144,77){\makebox(0,0){$y_0$}}
  \multiput(19,34)(0,15){2}{\line(1,0){6}}
  \multiput(144,34)(0,15){2}{\line(1,0){6}}
  \multiput(25,34)(0,15){2}{\circle*{2.8}}
  \multiput(144,34)(0,15){2}{\circle*{2.8}}
  \put(15,49){\makebox(0,0)[r]{$x_1$}}\put(154,49){\makebox(0,0)[l]{$y_1$}}
  \put(15,34){\makebox(0,0)[r]{$x_2$}}\put(154,34){\makebox(0,0)[l]{$y_2$}}%
}
\begin{picture}(358,83)
  \put(0,0){\begin{picture}(169,83)
    \cutpanel
    \thicklines
    \qbezier(25,71)(52,58)(85,60)
    \qbezier(85,60)(116,62)(144,71)
    \put(85,66){\makebox(0,0){$P$}}
    \thinlines
    \put(25,49){\line(1,0){119}}\put(25,34){\line(1,0){119}}
    \put(85,6){\makebox(0,0){old society}}
  \end{picture}}
  \put(189,0){\begin{picture}(169,83)
    \cutpanel
    \thicklines
    \qbezier(25,71)(54,58)(85,15)
    \put(85,15){\circle*{3}}
    \put(85,6){\makebox(0,0){$w_t$}}
  \end{picture}}
\end{picture}
\caption{The annulus opened along a radial arc, shown for $b=3$.
Its vertical sides are the two copies of the cut. The old vortex
attaches along the bottom society and is omitted from the drawing.
On the left, $P$ joins $x_0$ to $y_0$ above the other two paths of
an $X$--$Y$ linkage. On the right, a hypothetical first contact of
the path from $x_0$ with the society at $w_t$ gives a crosscut.
The remaining $X$ stubs lie to its left and the remaining $Y$ stubs
to its right. A path joining these sides through the vortex must meet
the separator at $w_t$.}
\label{fig:opened-vortex-annulus}
\end{figure}

Cut the annulus along $L$. The $b$ cycles minus their respective
contacts are disjoint paths. Denote their ends by $x_i,y_i$ for
$0\le i<b$, in order from the outermost cycle inwards. The portions of the cycle edges joining these ends to the contact give
disjoint stubs on opposite copies of the cut, on the $x$-side and
$y$-side, respectively. Each cycle has at least three vertices in the simple graph $F'$, so its two
ends are distinct and both survive, since $L$ meets the cycle only at its deleted contact. Let $X\coloneqq \{x_i:0\le i<b\}$
and $Y\coloneqq \{y_i:0\le i<b\}$. Consider the part drawn inside the outermost
cycle together with its old attachment, after the indicated deletions.
Every family of $b$ vertex-disjoint $X$--$Y$ paths in this graph has an
$x_0$--$y_0$ path entirely in the drawing and disjoint from the society.
Indeed, if the path from $x_0$ first hits the society at $w_t$, extend its
initial segment by the short $x_0$ stub. This gives a crosscut from the
outer corner of the opened annulus to $w_t$. Earlier contacts with the
outer boundary can be moved slightly into the annulus in sectors avoiding
the other paths. The crosscut separates the remaining $X$ stubs from the
$Y$ stubs. The two portions of the society on opposite sides of $w_t$ are
the earlier and later intervals of the opened society. Each of the other $b-1$ paths
must pass between those intervals through the old attachment and hence meets the
separator. The path from $x_0$
itself uses $w_t$, contradicting the separator bound $2a+1<b$. Finally, if
this drawing path ends at $y_i$ with $i>0$, the path with its two end stubs
cuts off $y_0$ from every remaining $X$ stub, with the entire old society
on their side. This contradicts the existence of the path to $y_0$.
Closing the $x_0$--$y_0$ path by the two outer contact edges bounds a disk
containing the other paths and the old attachment.
Figure~\ref{fig:opened-vortex-annulus} illustrates the crosscut.

The punctured cycles give an initial full $X$--$Y$ linkage. For each
internal vertex $v$ of its outer path, use a unit vertex-capacity flow
to test whether a full linkage exists after deleting $v$. Use only the
drawing in the current closed disk and its old attachment. If a test
succeeds, the new outer path bounds a smaller disk by the crosscut
argument. This disk omits $v$ since its boundary avoids $v$ and it lies
in the old disk. Thus at most $|V(G)|$ successful replacements are possible.

Write $J^\circ$ for the last disk graph including its attachment, and
$P\coloneqq p_0\cdots p_\ell$ for its outer path. Deleting any internal $p_i$ destroys the full linkage. Menger's theorem
gives a separator of size at most $b-1$ in $J^\circ-p_i$. Adding $p_i$
gives a minimum $X$--$Y$ separator of size $b$, since $J^\circ$ has a
linkage of size $b$. For $0<i<\ell$ put
\[
 D_i\coloneqq X\cup\{p_0,\ldots,p_i\},\qquad
 E_i\coloneqq Y\cup\{p_i,\ldots,p_\ell\}.
\]

Every order-$b$ separator through $p_i$ meets each linkage path exactly once.
The prefix of $P$ before $p_i$ is therefore on the $X$ side, its suffix on the
$Y$ side. Thus $D_i,E_i$ admit a separation of order $b$. Every minimum
separation has this order and its separator contains their common vertex $p_i$.

At each position $i$, choose a minimum separation $(A_i,B_i)$ with
$D_i\subseteq A_i$, $E_i\subseteq B_i$ and $|A_i|$ as small as possible.
These cuts are nested. For $i<j$, their lower corner
$(A_i\cap A_j,B_i\cup B_j)$ is feasible at $i$, and their upper corner
$(A_i\cup A_j,B_i\cap B_j)$ is feasible at $j$. Their orders sum to $2b$,
so both have order $b$. Minimality of $|A_i|$ forces $A_i\subseteq A_j$.
Since the lower corner has the same left side and order as the cut at $i$,
their right sides also agree, giving $B_j\subseteq B_i$.
To find the cut at $i$, enumerate sets of $b$ vertices as separators.
Discard a set if a surviving component meets both $D_i$ and $E_i$.
Otherwise put the separator on both sides and only the components meeting
$D_i$ on the left, and minimize the resulting left size. Write
$T_i\coloneqq A_i\cap B_i$. At the ends append $(X,V(J^\circ))$ and
$(V(J^\circ),Y)$, with separators $T_0\coloneqq X$ and $T_\ell\coloneqq Y$.

Take as the new society all vertices of $P$ in its path order. The vertex sets $B_i\cap
A_{i+1}$ give its consecutive pieces with root sets $T_i,T_{i+1}$. Assign linkage segments
between consecutive cuts to these pieces and assign every remaining edge to a piece containing
both its ends. Such a piece exists by nesting and because no edge crosses a separation. A
vertex shared by two pieces lies in every intervening root set. The only society vertices in a
piece are its two end roots on $P$. The $b$ path segments provide all the required local
linkages, so this is a rooted presentation. Move the boundary of the new disk slightly outwards, off the
open edges of $P$, and along the two closing edges at the removed contact. Thus the edges
of $P$ lie in the new disk and belong to the attachment. It meets the
remaining drawing just at the society vertices, in order on one boundary arc. All drawing
edges outside this disk remain in the drawing.

To bound the new attachments' treewidth, let $Z$ consist of the apices
and all opening vertices. Each arc contributes only its society contact and its
$b$ cycle contacts.
For a piece $J_i$ of a new attachment put
\[
 I\coloneqq V(J_i)\setminus(T_i\cup T_{i+1}),\qquad
 S_i\coloneqq Z\cup T_i\cup T_{i+1}.
\]
In the intact drawing $F^+$ with the remaining vortices, deleting $S_i$
separates $I$ from all vertices outside the piece. Indeed other pieces
meet $J_i$ only at its roots, which include its two society vertices.

This separation lifts to $G$. The ports of every small piece form a
clique in $F^+$, so its surviving ports lie on at most one strict side.
Put the piece's private vertices on that side, or on either side if all
ports were deleted. Put each set-aside vortex on the exterior side,
where all its society vertices lie. The more than $s$ original vertices
counted outside the selected disks are still on this strict side. The root
sets lie in the attachments and any drawn opening vertex lies on an
opening arc or in its vortex society. As $|S_i|\le A_\beta+u(b+1)+2b\le K$,
unbreakability gives $|I|\le s$. The pieces in order form a path
decomposition, each bag of size at most $s+2b$, so the treewidth is at most
$s+2b-1$.

Enumerate visits to the opening vertices by Lemma~\ref{near:apex} on the
intact graph with the exceptional ports reserved. This gives patterns of
bounded order $k'$ in the graph with the opening vertices removed. With no
remaining vortex the instance is already drawn. Otherwise remove the interiors
of the new attachment disks from the surface. The resulting compact surface
$S$ has one boundary circle per attachment, with the vertices of its outer
path in order on a proper subarc. The paths have distinct ends. The attachments are pairwise disjoint. Indeed the last
disks of distinct vortices are disjoint, and distinct vortices have no common vertex,
since they meet only in the drawing and their societies lie in their own disks. The
attachments have depth $b$ and treewidth at most $s+2b-1$.
Lemma~\ref{near:vortices} gives queries of order at most
$\thmbound{near:vortices}{S,\Gamma,u,b,k'}$ in subdrawings on the capped
surface, of Euler genus at most $g$. Up to homeomorphism there are finitely many
possibilities for $S$, since it has Euler genus at most $g$ and at most $u$ boundary
circles, and we take the largest of these orders. All the enumerations and flow computations
are polynomial for the fixed parameters. This completes the reduction.

\end{document}